\documentclass[aos, preprint]{imsart}

\RequirePackage{amsthm,amsmath,amsfonts,amssymb}
\RequirePackage[numbers]{natbib}
\RequirePackage[colorlinks,citecolor=blue,urlcolor=blue]{hyperref}
\RequirePackage{graphicx}
\usepackage{tikz}
\usepackage{pgfplots}
\pgfplotsset{compat=1.16}
\usepackage{color}

\startlocaldefs

\newtheorem{theorem}{Theorem}[section]
\newtheorem{corollary}[theorem]{Corollary}

\newtheorem{lemma}[theorem]{Lemma}
\theoremstyle{remark}

\newcommand{\Pois}{\mathrm{Pois}}

\newcommand{\Bin}{\mathrm{Bin}}
\newcommand{\Unif}{\mathrm{Unif}}
\newcommand{\LmbW}{\mathrm{W}_{-1}}
\newcommand{\Scal}{\mathcal{S}}

\newcommand{\reals}{\mathbb R}
\newcommand{\simiid}{\overset{\mathsf{iid}}{\sim}}

\newcommand{\Ncal}{\mathcal{N}}

\newcommand{\dense}{\mathsf{high}}
\newcommand{\sparse}{\mathsf{low}}

\newcommand{\Bonf}{\mathsf{Bonf}}
\newcommand{\one}{\mathbf{1}}
\newcommand{\HC}{\mathrm{HC}}
\newcommand{\minP}{\min\mathrm{P} }

\newcommand{\onesample}{\mathsf{one-sample}}
\newcommand{\twosample}{\mathsf{two-sample}}

\newcommand{\ex}[1]{\ensuremath{\mathbb{E}\left[ #1\right]}}

\newcommand{\exsub}[2]{\ensuremath{\mathbb{E}_{#1}\left[ #2\right]}}

\newcommand{\BonfColor}{cyan!80!white}
\newcommand{\normalColor}{green!60!black}
\newcommand{\HCcolor}{red!90!black}
\newcommand{\HCtextcolor}{red!40!white}
\newcommand{\normalTwoSmpColor}{\HCcolor}
\newcommand{\chisqColor}{orange}

\newcommand{\newtext}[1]{{\textcolor{black}{#1}}}

\endlocaldefs

\begin{document}

\begin{frontmatter}

\title{Higher Criticism to Compare Two Large Frequency Tables, with Sensitivity to Possible Rare and Weak Differences}


\runtitle{Higher Criticism for Two-sample Testing of Large Frequency Tables}

\begin{aug}
\author[A]{\fnms{David L.} \snm{Donoho}\ead[label=e1]{donoho@stanford.edu} }
\and
\author[B]{\fnms{Alon} \snm{Kipnis}\ead[label=e2]{alon.kipnis@idc.ac.il} }

\address[A]{Department of Statistics,
Stanford University,
\printead{e1}}
\address[B]{School of Computer Science,
Reichman University,
\printead{e2}}
\end{aug}

\begin{abstract}
We adapt Higher Criticism (HC) to the comparison of two frequency tables which may -- or may not -- exhibit  moderate differences between the tables in some unknown, relatively
small subset out of a large number of categories.
Our analysis of the power of the proposed HC test quantifies
the rarity and size of assumed differences and applies moderate deviations-analysis to determine the asymptotic powerfulness/powerlessness of our proposed HC procedure.

Our analysis considers the null hypothesis of no difference in underlying generative model
against a rare/weak perturbation alternative, in which the frequencies of $N^{1-\beta}$ out of the $N$ categories are perturbed by $r(\log N)/2n$ in the Hellinger distance; here $n$ is the size of each sample. 
Our proposed Higher Criticism (HC) test for this setting uses P-values obtained from $N$ exact binomial tests. We characterize the asymptotic performance of the HC-based test in terms of the rarity parameter $\beta$ and the perturbation intensity parameter $r$. Specifically, we derive a region in the $(\beta,r)$-plane where the test asymptotically has maximal power, while having asymptotically no power outside this region. Our analysis distinguishes between cases in which the counts 
in both tables are low, versus cases in which counts are high, corresponding to the cases of sparse and dense frequency tables. The phase transition curve of HC in the high-counts regime matches formally the curve delivered by HC in a two-sample normal means model.
\end{abstract}



\end{frontmatter}

\section{Introduction}
\label{sec:intro}

\subsection{Comparing Two Large Tables of Counts}
Suppose we have two frequency tables, each indexed by the same (large) collection of  categories, and we wish to know whether the underlying generating mechanisms behind the two tables might be different. Such a difference might indicate changes in time or changes caused by an intervention.  In addition, suppose that, if the generating mechanisms are different, we suspect  differences will arise only in a small fraction of the observed categories, but we do not know in advance where the differences are likely to occur. 

This question arises in a number of stylized applications, including:
\begin{itemize}
\item \emph{Attributing authorship}. We have two text corpora; do they have different authors? Changes between authors usually occur in word-frequencies of certain author-specific words \cite{MostellerWallace,kipnis2019higher}; however, there may be no fixed ``giveaway'' words, i.e., specific words whose use frequency determines authorship. Instead, there may be a list of words used slightly differently by a specific author, and that list may vary idiosyncratically from author to author; hence we don’t know in advance which words are likely to be informative. 
\item \emph{Public health surveillance}. Syndromic surveillance for early detection of health crises relies on anomalous behavior of counts in tables which span a large number of categories based on, say, clinical case features \cite{bravata2004systematic,miller2004syndromic}. 
An emergent public health situation may affect primarily frequencies in some subset of categories, but this subset varies from episode to episode.
\end{itemize}

For such applications, we want a tool to compare the two frequency tables that is particularly adapted to identifying changes in generating mechanism affecting a relatively small fraction of {\it a priori} unspecified categories.

This paper develops such a tool -- adapting Higher Criticism \cite{donoho2015special}
to this two-tables-of-counts setting.
Using P-values deriving from a Binomial Allocation model, our proposal is interpretable and easy to apply; our theoretical analysis shows that it is powerful against the above-mentioned changes in generating mechanism.

\newtext{Several recent works considered two-sample problems under a general ``closeness testing'' formulation \cite{acharya2012competitive,chan2014optimal,diakonikolas2016new,diakonikolas2019collision}. As opposed to these works, we study the ability to test for differences between distributions when, if there are differences, they will be hiding among a relatively few randomly scattered categories and they will be only moderately strong. It is well understood by now through previous work in other settings that the tests that are optimal for general closeness testing will not be able to be powerful under these very specific departures \cite{janssen2000global}.}

\subsection{Analysis Framework}
\label{sec:analysis_framework}
Our analysis will be conducted in the following mathematical framework. We consider two samples, each obtained by \newtext{roughly} $n$ independent draws from two possibly different distributions over the same finite set of $N$ categories. We would like to test whether the two distributions are identical, or not. Consider a \emph{rare/weak perturbation alternative}, where the difference between the two distributions are largely concentrated to a small, but unknown, subset of the $N$ categories. 
Specifically, let 
$(X_i)_{i=1}^N$ and $(Y_i)_{i=1}^N$ be the observed counts. The null and alternative have the following structure:
\begin{align}
    \begin{split}
    H_0 &\quad : \quad X_i,Y_i \sim \Pois(nP_i),\quad i=1,\ldots,N.\\
    H_1 &\quad : \quad X_i \sim \Pois(nP_i) \quad \text{and} \\
    &Y_i \sim (1-\epsilon)\Pois(nP_i) + \frac{\epsilon}{2} \Pois(nQ_i^+) + \frac{\epsilon}{2} \Pois(nQ_i^-)
    \quad i=1,\ldots,N.
        \label{eq:hyp}
    \end{split}
\end{align}
Here $P = (P_1,\ldots,P_N) \in \reals_+^N$ is a vector of `baseline' rates such that $\sum_{i=1}^N P_i = 1$, and the perturbations $Q_i^\pm$ obey
\begin{align}
\label{eq:perturbation_Hellinger}
    \sqrt{Q_i^\pm} \equiv  \max\left\{\sqrt{P_i} \pm \sqrt{\mu},0\right\},
\end{align}
for some $\mu \geq0$ to be determined later. The mixing fraction $\epsilon>0$ is typically small, and $\mu$ is relatively small as well. Informally, the case we explore in this article chooses $\epsilon$ and $\mu$ so that no single category $i$ can provide decisive evidence against the null hypothesis of identical distributions; the evidence is rare and weak. \par
Our analysis takes place in an asymptotic setting where both the number of features $N$ as well as the sample size $n$ go to infinity, but perhaps at different rates. 
We choose $\epsilon$ and $\mu$ according to $N$ and $n$ with:
\begin{subequations}
\label{eq:calibration}
\begin{align}
\label{eq:eps_def}
 \epsilon & = \epsilon_N \equiv N^{-\beta}, \quad \beta \in (0,1),
 \end{align}
 and
\begin{align}
\mu & = \mu_{N,n} \equiv r  \frac{\log(N)}{2n}, \quad r > 0. \label{eq:mu_def}
\end{align}
\end{subequations}
The parameter $(\beta,r)$ defines a \emph{phase space} of different situations,
\begin{itemize}
    \item[-] $\beta$ controls the rarity of the perturbation to be detected, with severe rarity $\beta \in (1/2,1)$ of most interest to us. 
    \item[-] $r$ controls the amplitude or intensity of the perturbation; the logarithmic calibration makes the testing problem $H_0$ vs. $H_1$ challenging but not impossible.
\end{itemize}

\par Our analysis of the testing problem \eqref{eq:hyp} behaves very differently depending on whether the counts of categories are typically high or low. The Low (respectively, High) Counts Case corresponds to a situation where the contingency table describing the data is sparse (respectively, dense). 
We discuss the two regimes separately.

\subsubsection*{High Counts Case}
In a \emph{high counts} scenario, the following conditions hold:
\begin{itemize}
    \item[]($\dense$) $\qquad \forall M>0$, $i=1,\ldots,N$,\quad 
    $\Pr\left( \frac{ \newtext{\min\{X_i,Y_i \}} } {\log(N)} > M \right) \to 1$,
\end{itemize}
under either $H_0$ or $H_1$.
Under this situation, the validity of
($\dense$) is determined by the underlying vector of rates $(P_1,\ldots,P_N)$ and the sample size $n$, but is unaffected by the rare-weak model parameters $\beta$ and $r$. Therefore, in terms of the rare-weak perturbation model \eqref{eq:hyp}, ($\dense$) is equivalent to the condition
\begin{align}
    \label{eq:dense_condition_P}
\frac{nP_i}{\log(N)} \to \infty.
\end{align}
From \eqref{eq:dense_condition_P} we also have 
\begin{align}
\label{eq:chi_square_perturb}    
(nQ_i^\pm \pm nP_i)^2 = 2 r \log(N) \cdot nP_i + o(1).
\end{align}
Hence, the perturbation is globally proportional to $\sqrt{P_i}$. Perturbations of this kind are very natural in statistics, in view of the important role of the Chi-squared and Hellinger discrepancies. Indeed, the typical term in the $\chi^2$-discrepancy, $(Q_i^\pm \pm P_i)^2/P_i$,
would equal simply $4\mu_{N,n} = 2 r \log(N)/n$ under such a perturbation model; hence the perturbation is naturally controlled in a Chi-squared sense between two rate vectors. \par
We note that the two-sample testing problem is symmetric in the two samples: each distribution might be seen as a perturbed version of the other when the vector of frequencies $P$ is unknown. In this paper, we sometimes speak of the $Y$ counts as being associated with the perturbed distribution, however this is from the point of view of our theoretical study, not from the practical viewpoint. 
\par
 
\subsubsection*{Low Counts Case}
Under a \emph{low counts} scenario, we have:
\begin{itemize}
    \item []($\sparse$) \qquad $\Pr \left( \newtext{\max\{X_i,Y_i\}} \leq \log(N)  \right) \to 1$, \qquad $i=1,\ldots,N$,
\end{itemize}
under either $H_0$ or $H_1$. By Markov's inequality, for $j\in\{0,1\}$,
\[
\newtext{\Pr_{H_j^{(N)}}\left( \max\{X_i,Y_i\} > \log(N)  \right) } \leq \Pr_{H_j^{(N)}}\left( X_i+Y_i > \log(N)  \right) \leq \frac{\ex{X_i+Y_i}}{\log(N)} = \frac{nP_i}{\log(N)}\left(2 + o(1)\right),
\]
and hence ($\sparse$) holds whenever: 
 \begin{equation}
\label{eq:sparse_condition_P}
\frac{nP_i}{\log(N)} \to 0;
\end{equation}
under which case, we have:
\begin{align}
    \label{eq:sparse_perturbation}
nQ_i^+ = nP_i +  \frac{1}{2}r\log(N) \left(1+o(1)\right), \qquad Q_i^- = 0.
\end{align}
In the Low Counts Case, the problem is not symmetric in the two samples. Indeed, the sample $(Y_1,\ldots,Y_N)$ has approximately $N^{1-\beta}$ entries with Poisson rates exceeding $r\log(N)/2$, while the number of entries in $(X_1,\ldots,X_N)$ with Poisson rates larger than $r\log(N)/2$ is smaller than $N^{\delta}$ for any $\delta>0$. 
\\

To fix ideas, we provide two examples for the baseline rates $P$ and the conditions under which the problem belongs to $(\dense)$ or $(\sparse)$:
\begin{itemize}
\item[(i)] {\bf Uniform baseline rates}: With $P_i=1/N$ and $n =N^{\xi}$, ($\dense$) holds if $\xi > 1$; ($\sparse$) holds if $\xi < 1$. 

\item[(ii)] {\bf Zipf-Mandelbrot baseline rates}: Assume that $P_i = c_N\cdot(i+k)^{-\xi}$ for some $k>-1$ and $\xi>1$; here $c_N$ is a normalization constant that satisfies
\[
c_N = \frac{1}{\sum_{i=1}^N (i+k)^{-\xi}}.
\]
$c_N$ is bounded away from zero since the sum in the denominator converges. We have
\[
\frac{n P_i}{\log(N)} = \frac{c_N}{\log(N)} \frac{n} {(i+k)^{\xi}} \geq \frac{c_N}{\log(N)} \frac{n} {(N+k)^{\xi}}
\]
and thus ($\dense$) holds if $n = N^{\gamma}$ and $ \xi <\gamma$, while ($\sparse$) holds if $\xi > \gamma$. 
\end{itemize}



\subsection{Binomial Allocation P-values}
The so-called `exact binomial test' P-value \cite{clopper1934use} is a function of $x$ and $y$, where for $x,y \in \mathbb N$ we set 
\begin{align}
\label{eq:pval_def_pre}
\pi(x,y) & \equiv 
 \Pr\left( \left|\Bin(n,p) - np\right| \leq \left|x-np \right| \right);
\end{align}
here $p=1/2$ and $n\equiv x+y$ (note the symmetry $\pi(x,y)=\pi(y,x)$ with this choice of $p$ and $n$). Our testing procedure based on $\pi_i$ or $\tilde{\pi}_i$, $i=1,\ldots,N$, can be easily modified to address cases where both samples have non-equal sizes; see \cite{kipnis2019higher,kipnis2021two} for the details. The (exact) P-value associated with the $i$-th feature (category) is
\begin{align}
\pi_i \equiv \pi(X_{i},Y_{i}). \label{eq:Pval_def}
\end{align}
We also consider the randomized P-value $\tilde{\pi}_i$ associated with $\pi_i$. This P-value is uniformly distributed under the null and is dominated by $\pi_i$ both under the null and alternative. Specifically,
\begin{align}
\label{eq:Pval_def_rand}
    \Pr_{H_0}\left(\tilde{\pi}_i \leq t \right) = t,\quad t\in [0,1], \quad \text{and}\quad \Pr_{H_j^{(N)}}\left( \tilde{\pi}_i \leq \pi_i \right) = 1,\quad j\in\{0,1\}.
\end{align}

\subsection{Higher Criticism}
We combine the collection of P-values $\pi_1,\ldots,\pi_N$ into a global test against $H_0$ by applying Higher Criticism \cite{donoho2015special}. Define the HC component score:
\[
\HC_{N,n,i} \equiv \sqrt{N} \frac{i/N - \pi_{(i)}}{\sqrt{\pi_{(i)}\left(1-\pi_{(i)}\right)}},
\]
where $\pi_{(i)}$ is the $i$-th ordered P-value among $\{ \pi_i,\, i=1,\ldots,N \}$. The HC statistic is:
\begin{align}
    \label{eq:HC}
\HC_{N,n}^\star \equiv 
\underset{{1\leq i \leq N\gamma_0}}{\max}  \HC_{N,n,i},
\end{align}
where $0<\gamma_0 < 1$ is a tunable parameter\footnote{ $\gamma_0$ typically has no effect on the asymptotic value of $\HC_{N,n}^\star$ under $H_1$. Often $\gamma_0 = 1/20$ or $\gamma_0 = 1/10$.
}. \par
We reject $H_0$ at level $\alpha$ when $\HC_{N,n}^\star$ exceeds the $95$ percentile (say) of under the null. 

\subsection{Performance of HC Test}

The power of the test varies dramatically across the $(\beta,r)$ phase space. In part of this region, the test will work well; in another part it will fail to detect. Formally, for a given sequence of statistics $\{T_{N,n}\}$ and hypothesis testing problems \eqref{eq:hyp} indexed by $n$ and $N$ where $N=N(n)$, we say that $\{T_{N,n}\}$ is \emph{asymptotically powerful} if there exists a sequence of thresholds $\{h(n)\}$ such that
\[
\Pr_{H_0} \left( T_{N,n} > h(n) \right) + \Pr_{H_1} \left( T_{N,n} \leq h(n) \right) \to 0,
\]
as $n$ goes to infinity. In contrast, we say that $\{T_{N,n}\}$ is \emph{asymptotically powerless} if 
\[
\Pr_{H_0} \left( T_{N,n} > h(n) \right) + \Pr_{H_1} \left( T_{N,n} \leq h(n) \right) \to 1,
\]
for any sequence $\{h(n)\}_{n\in \mathbb N}$. \par
The statistic $\HC_{N,n}^\star$ experiences a \emph{phase transition} in $(\beta,r)$: for a specific function $\rho(\beta)$ given below, 
$\HC_{N,n}^\star$ is asymptotically powerful when $r > \rho(\beta)$ and asymptotically powerless when $r < \rho(\beta)$. Our main results characterize the function $\rho(\beta)$ under each of the cases ($\dense$) and ($\sparse$).

\subsubsection{High Counts Case}
Define the would-be phase transition boundary
\begin{align}
\label{eq:rho_dense_def}
\rho_{\dense}(\beta) \equiv \begin{cases} 2(\beta-1/2) &  1/2\leq \beta < 3/4, \\
2(1 - \sqrt{1-\beta})^2 &  3/4 \leq \beta \leq 1.
\end{cases}
\end{align}

Figure~\ref{fig:phase_diagram_dense} illustrates the curve $\rho_{\dense}(\beta)$. 
\begin{figure}
    \begin{center}
    \begin{tikzpicture}
    \begin{axis}[
    width=10cm,
    height=7cm,
    legend style={at={(0,1)},
      anchor=north west, legend columns=1},
    ylabel={$r$ (intensity)},
    xlabel={$\beta$ (rarity)},
    ytick={0,0.17,0.5,1,1.5,2},
    yticklabels={0,\scriptsize $3-2\sqrt{2}$,0.5,1,1.5,2},
    xtick={0.5,0.6,0.7,0.75,0.8,0.9,1},
    xticklabels={.5,.6,.7,.75,.8, .9,1},
    legend cell align={left},
    ymin=0,
    xmin=0.5,
    xmax=1,
    ymax=2,
    ]

\addplot[domain=0.4:1, color=\BonfColor, samples=121, style=ultra thick] {2*(sqrt(1)-sqrt(1-x))^2};
\addlegendentry{$\rho_{\dense}^{\Bonf}(\beta)$ ($\minP$ test)}

\addplot[domain=0.75:1, samples = 121, color=\normalTwoSmpColor, style=ultra thick]
    {2*(1-sqrt(1-x))^2};
    \addlegendentry{$\rho_{\dense}(\beta)$ (HC test)};

\addplot[domain=0.5:0.75, color=\normalTwoSmpColor, style=ultra thick] {2*(x-1/2)};

\addplot[dotted, style=thick] 
coordinates {(0.75,0) (0.75,0.5)};

\addplot[dotted, style=thick] 
coordinates {(0,0.5) (0.75,0.5)};




\node (succ) at (axis cs:.8,1){};
\node (fails) at (axis cs:.85,.3){};
\end{axis}

\node[fill = \HCtextcolor, rotate=15] at (succ) {$\HC$ succeeds};

\node[fill = \HCtextcolor, rotate = 15] at (fails) {$\HC$ fails};

\end{tikzpicture}
\end{center}
    
\caption{Phase Diagram (High Counts Case):
The phase transition curve $\rho_{\dense}(\beta)$ of \eqref{eq:rho_sparse_def} separates between the region where HC test is asymptotically powerful and asymptotically powerless. Also shown is the region of success for the min-P test used in Bonferroni-type inference.
} 
\label{fig:phase_diagram_dense}
\end{figure}


\begin{theorem}[High Counts Case] \label{thm:HC_power1}
Consider Problem \eqref{eq:hyp} under ($\dense$) with parameters $\beta$ and $\epsilon$ calibrated with $n$ and $N$ as in \eqref{eq:eps_def} and \eqref{eq:mu_def}. 
\begin{itemize}
    \item[(i)] {\it Region of Full Power}: Assume that $r> \rho_{\dense}(\beta)$. Consider the Higher Criticism statistic based either on randomized P-values \eqref{eq:Pval_def_rand} or non-randomized P-values \eqref{eq:Pval_def}. Then through this region, $\HC_{N,n}^\star$ is asymptotically powerful. 
    \item[(ii)] {\it Region of No Power}: Instead assume that $r <  \rho_{\dense}(\beta)$ and consider the Higher Criticism statistic based on the randomized P-values \eqref{eq:Pval_def_rand}. Then throughout this region, $\HC_{N,n}^\star$ is asymptotically powerless.
\end{itemize}
\end{theorem}

\begin{figure}
    \begin{tikzpicture}
    \begin{axis}[
    width=10cm,
    height=7cm,
    legend style={at={(0,1)},
      anchor=north west, legend columns=1},
    ylabel={$r$ (intensity)},
    xlabel={$\beta$ (rarity)},
    xtick={0.4,0.5,0.6,0.7,0.8,0.9,0.923,1},
    ytick={0,1,2,2.04,3},
    yticklabels={0,1,,$\frac{\sqrt{2}}{\log(2)}$,3},
    xticklabels={.4,.5,.6,.7,.8,,$\frac{1}{2} + \frac{1-\frac{1}{\sqrt{2}}}{\log(2)}$,1},
    legend cell align={left},
    ymin=0,
    xmin=0.5,
    xmax=1,
    ymax=3]

\addplot[domain=1:2.89, color=\BonfColor, samples=121, style=ultra thick]
    ({1 + (2 - x*ln(2) - 2* ln(2/(x * ln(2))))/ln(4)},{x});
\addlegendentry{$\rho_{\sparse}^{\Bonf}$ ($\minP$ test)}

\addplot[domain=0.5:0.923, color=\HCcolor, style=ultra thick] 
    {4.83*(x-1/2)};
    \addlegendentry{$\rho_{\sparse}$ ($\HC$ test)}

\addplot[domain=2:2.89, color=\HCcolor, samples=121, style=ultra thick]
    ({1 + (2 - x*ln(2) - 2* ln(2/(x * ln(2))))/ln(4)},{x});

\addplot[color=black, style=thick,dotted]
coordinates {(1,2.88539) (0.98,2.88539)}; 

\addplot[color=black, style=thick,dotted]
coordinates {(0.923,0) (0.923,2.045)}; 

\addplot[color=black, style=thick,dotted]
coordinates {(0,2.045) (0.923,2.045)}; 



\node[] (topl) at (axis cs: 0.4,2.04) {};
\node[] (topr) at (axis cs: 1,2.885) {};
\node[] (bottomr) at (axis cs: 0.923,0.2) {};

\end{axis}


\node[xshift = 0.7cm] at (topr) {$\frac{2}{\log(2)}$};

\node[fill = \HCtextcolor, rotate = 28] at (2,1.5) {$\HC$ succeeds};

\node[fill = \HCtextcolor, rotate = 28] at (2.3,0.7) {$\HC$ fails};

 
\end{tikzpicture}
\caption{Phase Diagram (Low Counts Case):
The phase transition curve $\rho_{\sparse}(\beta)$ of \eqref{eq:rho_sparse_def} separates between the region where the HC test is asymptotically powerful and asymptotically powerless. Also shown is the region of success for the min-P test used in Bonferroni-type inference.
}
\label{fig:phase_diagram_sparse}
\end{figure}

\subsubsection{Low Counts Case}

Define 
\begin{align}
    \label{eq:rho_sparse_def}
    \rho_{\sparse}(\beta) \equiv \begin{cases}
    2 (1 + \sqrt{2})\left(\beta-\frac{1}{2}\right) & \frac{1}{2} < \beta \leq \frac{1}{2}+\frac{\sqrt{2}-1}{\sqrt{2}\log (2)}, \\
    -\frac{2 \LmbW\left(-\frac{2^{-\beta}}{2e}\right)}{\log (2)}&  \frac{1}{2}+\frac{\sqrt{2}-1}{\sqrt{2}\log (2)} < \beta < 1,
    \end{cases}
\end{align}
where $\LmbW(x)$ is the real negative branch of the Lambert $W$ function\footnote{In an earlier version of the manuscript our expression for $\rho_{\sparse}(\beta)$ was not written in terms $\LmbW$. We thank one of the referees for bringing to our attention that this expression can be written using the Lambert function. }, obtained as the negative solution $y$ of $x=ye^y$ \cite{corless1996lambertw}, \cite[Exc. 2.18]{LugosiConcentration}.

\begin{theorem}[Low Counts Case] \label{thm:HC_power2}
Consider Problem \eqref{eq:hyp} under ($\sparse$) with parameters $\beta$ and $\epsilon$ calibrated with $n$ and $N$ as in \eqref{eq:eps_def} and \eqref{eq:mu_def}. 
\begin{itemize}
    \item[(i)] Region of Full Power: Assume that $r> \rho_{\sparse}(\beta)$. Consider the Higher Criticism statistic based either on randomized P-values \eqref{eq:Pval_def_rand} or non-randomized P-values \eqref{eq:Pval_def}. Then through this region, $\HC_{N,n}^\star$ is asymptotically powerful. 
    \item[(ii)] Region of No Power: Instead assume that $r <  \rho_{\sparse}(\beta)$ and consider the Higher Criticism statistic based on the randomized P-values \eqref{eq:Pval_def_rand}. Then throughout this region, $\HC_{N,n}^\star$ is asymptotically powerless.
\end{itemize}
\end{theorem}

\subsection{Optimality of the Phase Diagram}
Theorems~\ref{thm:HC_power1} and \ref{thm:HC_power2} provide precise descriptions of the phase diagram of HC in the two-sample setting \eqref{eq:hyp}. These theorems also suggest that HC is asymptotically powerful for any $r>0$ when $\beta < 1/2$, a fact that easily follows from the discussion in \cite[Sec. 6.1]{cai2014optimal}.

Regarding the broader question of the optimality of HC's phase diagram, there is a larger context that is developed in \cite{kipnis2021logchisquared} showing that  $\rho_{\dense}(\beta)$ describes the optimal phase transition for the High Counts Case. Namely, any test statistic in the High Counts Case that is based on the P-values \eqref{eq:Pval_def_rand} is asymptotically powerless for $r < \rho_{\dense}(\beta)$. No equivalent result is known in the Low Counts Case. 

\subsection{The one-sample Poisson rare/weak setting}
\label{sec:one_sample}
Arias-Castro and Wang \cite{arias2015sparse} studied the goodness-of-fit problem of the Poisson rates $(\lambda_1,\ldots,\lambda_N)$ and the sample $(Y_1,\ldots,Y_N)$, where:
\begin{align} 
\begin{split}
\label{eq:one_sample}
    H_0 \, : & Y_i \simiid \Pois(\lambda_i), \quad i=1,\ldots,N,  \\
    H_1 \, : & Y_{i} \simiid (1-\epsilon_n)\Pois(\lambda_i) + \frac{\epsilon_n}{2} \Pois(\lambda_i^+) + \frac{\epsilon_n}{2} \Pois(\lambda_i^-), \\
    & \quad i=1,\ldots,N.
    \end{split}
\end{align}
They considered two different regimes for the parameters $\lambda_i$, $\lambda_i^+$ and $\lambda_i^-$:
\begin{itemize}
    \item \emph{Large} Poisson means: 
    \[
    \lambda_i / \log(N) \to \infty,  \quad \text{and} \quad \lambda_i^\pm = \lambda_i \pm \sqrt{2r\log(N) \lambda_i}.
    \]
    \item \emph{Small} Poisson means:
    \[
    \lambda_i / \log(N) \to 0, \quad \lambda_i^+ = \lambda_i^{1-\gamma} (\log(N))^\gamma, \quad \text{ and } \lambda_i^- = 0.
    \]
\end{itemize}
Comparing with \eqref{eq:dense_condition_P}-\eqref{eq:chi_square_perturb} and \eqref{eq:sparse_condition_P}-\eqref{eq:sparse_perturbation}, we see that our distinction made here, between low and high counts 
, is analogous to the distinction between large and small Poisson means made in \cite{arias2015sparse}. It follows that, if the vector $P = (P_1,\ldots,P_N)$ underlying $H_0$ in \eqref{eq:hyp} is fully known to us, and if the data $(X_i)$ are unobserved by us, we obtain a modified testing problem that is essentially \eqref{eq:one_sample}.
We summarize the results of \cite{arias2015sparse} relevant to our setting. \par

\subsubsection{High sample case of the one-sample problem}
Define 
\begin{equation}
    \label{eq:rho_onesample}
\rho_{\onesample}(\beta) \equiv \begin{cases}
(1-\sqrt{1-\beta})^2 &  3/4 < \beta < 1, \\
\beta - \frac{1}{2} & 1/2 < \beta \leq 3/4.
\end{cases}
\end{equation}
The results of \cite{arias2015sparse} imply that $\rho_{\onesample}(\beta)$ describes a fundamental phase-transition for \eqref{eq:one_sample} under \eqref{eq:dense_condition_P}: all tests are asymptotically powerless when $r < \rho_{\onesample}(\beta)$, while some tests are asymptotically powerful whenever $r>\rho_{\onesample}(\beta)$. Specifically, a version of HC that uses P-values obtained from a normal approximation to the Poisson random variables is asymptotically powerful whenever $r > \rho_{\onesample}(\beta)$. \par
Note that 
\[
\rho_{\dense}(\beta) = 2\rho_{\onesample}(\beta),
\]
hence the two-sample phase transition for HC is at a different location than the one-sample phase transition. 
%

%
\subsubsection{Low Counts Case of the one-sample problem}
Under Case $(\sparse)$ \eqref{eq:sparse_condition_P}, it follows from \cite[Prop. 7 ]{arias2015sparse} that in the one-sample setting the min-P test, and hence also HC, achieves maximal power over the entire range\footnote{The setting of \cite{arias2015sparse} only considered the case $r=2$, but the proof there extends in a straightforward manner to any $r>0$.} $0<r$ and $0<\beta<1$. In view of Theorem~\ref{thm:HC_power2}, the distinction between the one-sample and the two-sample setting is much more dramatic in the Low Counts Case than the High Counts Case; HC in the two-sample setting has a non-trivial phase transition, while in the one-sample setting HC is asymptotically powerful over the entire phase plane $\{1/2 < \beta < 1,\,0<r<1\}$.

\subsection{The one-sample normal means model}
The work of Donoho and Jin \cite{donoho2004higher} studied the behavior of HC under the one-sample rare/weak normal means setting:
\begin{align}
\begin{split}
 H_0 & : \, Y_i \simiid \Ncal(0,1),\quad i=1,\ldots,N, \\   
 H_1 & : \, Y_i \simiid (1-\epsilon_N)\Ncal(0,1) + \epsilon_N \Ncal(\mu_N,1), \quad i=1,\ldots,N,
 \label{eq:one_sample_normal_means}
 \end{split} 
\end{align}
where $\epsilon_N = N^{-\beta}$ and $\mu_N = \sqrt{2 r \log(N)}$. Specifically, it was shown in \cite{donoho2004higher} that HC is asymptotically powerful within the entire range of parameters $(\beta,r)$ under which the problem \eqref{eq:one_sample_normal_means} is solvable. Articles \cite{ingster1996some} and \cite{jin2003detecting} derived this range to be  $r > \rho_{\onesample}(\beta)$ of \eqref{eq:rho_onesample}. Several studies of HC behavior in rare/weak settings analogous to \eqref{eq:one_sample_normal_means} also experience phase transitions described by $\rho_{\onesample}(\beta)$ \cite{jin2003detecting,ingster2010detection,tony2011optimal,arias2011global,mukherjee2015hypothesis,arias2015sparse}. \par
Slightly more relevant to our discussion than \eqref{eq:one_sample_normal_means} is a one-sample normal means model with a two-sided perturbation alternative:
\begin{align}
\begin{split}
 H_0 & : \, Y_i \simiid \Ncal(0,1),\quad i=1,\ldots,N, \\   
 H_1 & : \, Y_i \simiid (1-\epsilon_N)\Ncal(0,1) + \frac{\epsilon_N}{2}\Ncal(\mu_N,1) + \frac{\epsilon_N}{2}\Ncal(-\mu_N,1) \quad i=1,\ldots,N.
 \label{eq:one_sample_normal_means1}
 \end{split} 
\end{align}
Arguing as in \cite{donoho2004higher}, it is straightforward to verify that HC of the P-values 
\[
\tilde{\pi}_{i} = \Pr \left( \Ncal(0,1) \geq |Y_i| \right),\quad i=1,\ldots,N,
\]
has the phase transition given by $\rho_{\onesample}(\beta)$. Namely, HC has the same phase transition in problems \eqref{eq:one_sample_normal_means} and \eqref{eq:one_sample_normal_means1}. 

\subsection{The two-sample normal means model}

Consider a two-sample normal means model:
\begin{align}
\begin{split}
 H_0 & : \, X_i,Y_i \simiid \Ncal(\nu_i,1),\quad i=1,\ldots,N, \\
 H_1 & : \, X_i \simiid \Ncal(\nu_i,1),\quad \text{and}\label{eq:two_sample_normal_means}
 \\
 & Y_i \simiid (1-\epsilon_N)\Ncal(\nu_i,1) + \frac{\epsilon_N}{2} \Ncal(\nu_i^+ ,1) + \frac{\epsilon_N}{2} \Ncal(\nu_i^-,1),
 \quad i=1,\ldots,N,
 \end{split}
\end{align}
with the perturbations 
\begin{align}
    \label{eq:two_sample_normal_means_perturbation}
    \nu_i^{\pm}-\nu_i = \pm\sqrt{2 r \log(N)}.
\end{align} 
We will show that HC of the P-values
\begin{align}
    \label{eq:normal_Pvalues}
\bar{\pi}_i \equiv \Pr\left( |\Ncal(0,1)| \geq \frac{\left|Y_i-X_i \right|}{\sqrt{2}} \right),\quad i=1,\ldots,N,
\end{align}
has the same phase transition in the model \eqref{eq:two_sample_normal_means} as it has in \eqref{eq:hyp} under the High Counts Case. 
Note that, in analogy with the binomial P-values \eqref{eq:Pval_def}, $\bar{\pi}_1,\ldots,\bar{\pi}_N$ are obtained from the data without specifying the means $\nu_1,\ldots,\nu_N$; these means remain unknown to us\footnote{If $\nu_1,\ldots,\nu_N$ are known, subtracting them from $(Y_i)$ leads to the one-sample problem \eqref{eq:one_sample_normal_means}.}.
\begin{theorem} \label{thm:two_sample_normal_means}
Consider the two-sample problem \eqref{eq:two_sample_normal_means} where $\beta$ and $r$ are calibrated to $N$ as in \eqref{eq:eps_def} and \eqref{eq:mu_def}. The higher criticism of the P-values $\bar{\pi}_1,\ldots,\bar{\pi}_{N}$ is asymptotically powerful whenever $r > \rho_{\dense}(\beta)$ and asymptotically powerless whenever $r < \rho_{\dense}(\beta)$.
\end{theorem}
Consequently, 
\begin{corollary}
The asymptotic phase transition of higher criticism of the P-values \eqref{eq:normal_Pvalues} in the two-sample normal means problem \eqref{eq:two_sample_normal_means} is 
$2 \rho_{\onesample}(\beta)$. 
\end{corollary}

\subsection{Bonferroni/Min-P Test}

Like HC, Bonferroni inference uses all the P-values; however, it only explicitly uses the smallest P-value $\pi_{(1)}$. The following theorems derive the region where the \emph{min-P test}, i.e., a test relying on $\pi_{(1)}$, is asymptotically powerful.
\begin{theorem} \label{thm:Bonferroni_dense}
Define
\[
\rho_{\dense}^{\Bonf}(\beta) \equiv 2\left(1-\sqrt{1-\beta }\right)^2, \qquad 1/2 \leq \beta \leq 1.
\]
Consider the hypothesis setting \eqref{eq:hyp} with the binomial P-values $\pi_1,\ldots,\pi_N$ of \eqref{eq:Pval_def}. A test based on $\pi_{(1)} = \min_i \pi_{i}$ is asymptotically powerful whenever 
\[
r > \rho_{\dense}^{\Bonf}(\beta).
\]
\end{theorem}

\begin{theorem} 
\label{thm:Bonferroni_sparse}
Define
\[
\rho_{\sparse}^{\Bonf}(\beta) \equiv 
-\frac{2 \LmbW \left(-\frac{2^{-\beta}}{2e}\right)}{\log (2)},
 \qquad 1/2 \leq \beta \leq 1.
\]
Consider the hypothesis setting \eqref{eq:hyp} with the binomial P-values $\pi_1,\ldots,\pi_N$ of \eqref{eq:Pval_def}. A test based on $\pi_{(1)} = \min_i \pi_{i}$ is asymptotically powerful whenever 
\[
r > \rho_{\sparse}^{\Bonf}(\beta).
\]
\end{theorem}
As Figures~\ref{fig:phase_diagram_dense} and \ref{fig:phase_diagram_sparse} show, the min-P test has phase diagram equally good to HC on the segment $\beta > \beta_0$, where 
\[
\beta_0 = \begin{cases} 3/4 & \text{($\dense$)}, \\
\frac{1}{2} + \frac{1-\frac{1}{\sqrt{2}}}{\log(2)} & 
\text{($\sparse$)}.
\end{cases}
\]
Hence, under sufficient rarity, Bonferroni inference is just as good as HC.

\subsection{Structure of this paper}
Section~\ref{sec:analysis} below presents an heuristic discussion of Theorems~\ref{thm:HC_power1} and \ref{thm:HC_power2}. 
In Section~\ref{sec:sim}, we provide simulations to support our theoretical findings. Discussion and concluding remarks are provided in Section~\ref{sec:discussion}. 
All proofs are provided in Section~\ref{sec:main_proof}.




\section{Where does HC find the Evidence? \label{sec:analysis}}

Previous studies observed that HC implicitly identifies a specific, data-driven subset of the observed P-values as driving the decision to possibly reject $H_0$. Donoho and Jin observed, in a different setting, that this subset may serve as an optimal set of discriminating features \cite{donoho2008higher,donoho2009feature}. The location of the specific informative P-values varies with the model parameters $\beta$ and $r$ (HC is adaptive since these parameters need not be specified by us). \par
In the High Counts Case, the behavior of the binomial P-values is analogous to the normal P-values in the one-sample normal means model \eqref{eq:one_sample}, as discussed in \cite{donoho2004higher} and \cite{gontscharuk2015intermediates}. Their behavior is different in the Low Counts Case. \\

Consider two versions of the empirical CDF of the binomial P-values
\begin{align*}
F^-_{N,n}(t) \equiv \frac{1}{N}\sum_{i=1}^N \one\{\pi_i < t\},\qquad F_{N,n}(t) \equiv \frac{1}{N}\sum_{i=1}^N \one \{\pi_i \leq t\} 
\end{align*}
Of course these are the same except at jumps, where they are left-continuous and right-continuous, respectively.
Note that, for $\gamma_0 < 1/2$,
\begin{align}
    \label{eq:HC_sandwitch}
\max_{1/N\leq t \leq \gamma_0} \sqrt{N}\frac{F^-_{N,n}(t)-t}{\sqrt{t(1-t)}} \leq \HC_{N,n}^\star \leq \max_{\pi_{(1)} \leq t \leq  \gamma_0} \sqrt{N}\frac{F_{N,n}(t)-t}{\sqrt{t(1-t)}}.
\end{align}
Evidence for the difference between $H_0$ and $H_1$ is to be sought among the smallest P-values; we anticipate that both sides of \eqref{eq:HC_sandwitch} attain their maximum at values of $t$ approaching zero. Since $F^-_{N,n}(t)= F_{N,n}(t)$ almost everywhere, we focus our attention on evaluating the HC component score $\HC_{N,n,i}$ for $i/N$ small. Equivalently, we consider
\[
V(t_n) \equiv \sqrt{N}\frac{F_{N,n}(t_n)-t_n}{\sqrt{t_n(1-t_n)}} \sim \sqrt{N}\frac{F_{N,n}(t_n)}{\sqrt{t_n}} - \sqrt{Nt_n}
\]
where $\{t_n\}$ is a sequence that goes to zero slowly enough as $n$ and $N$ go to infinity. In what follows, we use the sequence $t_n = N^{-q}$ where $q>0$ is a fixed exponent. Under $H_0$, the P-values have a distribution that is close to uniform -- $\Pr_{H_0 }(\pi_i\leq t)\sim t$; hence
\[
\Pr_{H_0} \left( \pi_i \leq N^{-q} \right) = N^{-q+o(1)},
\]
where the notation $o(1)$ represents a deterministic sequence tending to zero as $n$ and $N$ go to infinity. Evaluating $F_{N,n}(N^{-q})$ using the last display, it follows that $V(N^{-q})$ is bounded in probability under the null. The theoretical engine driving our main results is the following characterization of the P-values under $H_1$:
\[
\Pr_{H_1} \left( \pi_i \leq N^{-q} \right) =  N^{-\beta-\alpha_*(q,r)+o(1)} + N^{-q+o(1)},
\]
where 
\begin{align}
\label{eq:alpha_def}
\alpha_{*}(q,r) \equiv 
\begin{cases} (\sqrt{q}-\sqrt{r/2})^2 & * =\dense, \\
q\frac{\log \left(\frac{2q}{ r\log(2)}\right)-1}{\log (2)} + \frac{r}{2} & * = \sparse.
\end{cases}
\end{align}
In the dense case, this characterization is given by Lemma~\ref{lem:pi_q_dense}. That lemma uses a Chernoff bound argument to approximate the binomial test and in this sense relies on normal approximation for the model \eqref{eq:hyp}. In the sparse case, the behavior of the binomial P-values is given by Lemma~\ref{lem:pi_q_sparse}; it approximates the binomial test \eqref{eq:Pval_def} for values of $x$ close to zero. 
Altogether, these results lead to:
\[
\exsub{H_1}{V(N^{-q})} \sim N^{\frac{q+1}{2}-\beta-\alpha_*(q,r)} - N^{\frac{1-q}{2}},\quad N \to \infty,
\]
where $* \in \{\dense,\sparse\}$. Roughly speaking, the most informative part of the data corresponds to the location $t^\star=N^{-q^\star}$, where $q^\star$ maximizes the growth rate of $\ex{V(N^{-q})}$ under $H_1$. More explicitly, define
\begin{align}
    \label{eq:Xi_def}
\Xi_*(q,\beta,r) \equiv \frac{q+1}{2} - \beta - \alpha_{*}(q,r),\quad *\in\{\sparse,\dense\}.
\end{align}
The phase transition curve $\rho_*(\beta)$ is the boundary of the phase diagram region $\{(r,\beta)\,:\, \Xi_*^\star(r,\beta)>0\}$, where
\begin{align}
    \label{eq:max_q}
\Xi_*^\star(r,\beta) \equiv \max_{0 \leq q\leq 1} \Xi(q,\beta,r). 
\end{align}
Our reason for restricting $q$ to values at most one in \eqref{eq:max_q} is that, under $H_0$, essentially no P-values smaller than $cN^{-1}$ will occur, so there is no need to ``look further out'' than $q=1$. 
\par
The boundary of the phase diagram region $\{(r,\beta)\,:\,\Xi_*^\star(r,\beta)>0\}$ behaves differently for $*=\dense$ or $*=\sparse$; below we consider each case separately. 

\subsubsection{High Counts Case}
We have
\[
\Xi^\star_\dense(r,\beta) =  \begin{cases}
 1 - \beta - (1-\sqrt{r/2})^2 & r \geq 1/2, \\
  \frac{1+r}{2} - \beta & r < 1/2,
\end{cases}
\]
with $q^\star_{\dense}(r)$ attaining the maximum in \eqref{eq:max_q} is given by
\begin{align}
q^\star_{\dense}(r) = \begin{cases}
 1 & r \geq 1/2, \\
 2r & r < 1/2. 
\end{cases}
\label{eq:q_star_dense}
\end{align}
In short, for $r$ greater than $1/2$, evidence against $H_0$ is found at P-values of size $\asymp$ $N^{-1}$; i.e., in the very smallest P-values. In this region, the phase diagram for HC is equivalent to the phase diagram for the min-P; see Theorem~\ref{thm:Bonferroni_dense} below. The situation is different, however, for values of $r$ smaller than $1/2$. In such cases, the most informative part of the data is given by P-values $\asymp$ $N^{-2r}$. 

\subsubsection{Low Counts Case}
We have 
\begin{align}
\Xi_{\sparse}^\star(r,\beta) = \begin{cases}
 -\beta -\frac{r}{2}+\frac{\log (r)+1+\log (\log (2))}{\log (2)}
  & r \geq \frac{\sqrt{2}}{\log(2)}, \\
    \frac{\sqrt{2}-1}{2} r -\beta +\frac{1}{2} & r < \frac{\sqrt{2}}{\log(2)},
\end{cases}
\end{align}
with $q^\star_{\sparse}(r)$ attaining the maximum in \eqref{eq:max_q} is given by
\begin{align}
q^\star_{\sparse}(r) = \begin{cases}
 1 & r \geq \frac{\sqrt{2}}{\log(2)}, \\
 r \log(2)/\sqrt{2} & r < \frac{\sqrt{2}}{\log(2)}. \label{eq:q_star_sparse}
\end{cases}
\end{align}
In particular, the function $\LmbW(x)$ in \eqref{eq:rho_sparse_def} arises when solving $\Xi^\star_{\sparse}(r,\beta)=0$ for $r$. The two regimes for $r$ in  \eqref{eq:q_star_sparse} are analogous to the two regimes of $r$ in \eqref{eq:q_star_dense}. Superficially, for $r>\sqrt{2}/\log(2)$ the boundary of the region $\{\Xi^\star_{\sparse}(r,\beta) > 0\}$ is the same as the boundary of the region where the min-P test is powerful. In contrast, in the region $r<\sqrt{2}/\log(2)$, the most informative part of the data depends on $r$ and is given by P-values of size $\asymp$ $N^{-r \sqrt{2}/\log(2)}$. \par

\newcommand{\testSizeDense}{0.55}
\newcommand{\testSizeSparse}{0.6}
\newcommand{\sigLevel}{0.05}
\newcommand{\npow}{5}
\newcommand{\xiSparse}{0.8}
\newcommand{\xiDense}{1.4}

\section{Simulations}
\label{sec:sim}
We now discuss numerical experiments illustrating our theoretical results.\par
Our experiments involve Monte-Carlo simulations at each point  $(\beta,r)$ in a grid $I_r \times I_\beta$ covering the range $I_r \subset [0,3]$, $I_\beta \subset [0.45,1]$. Here $\beta$ and $r$ are as in \eqref{eq:eps_def} and \eqref{eq:mu_def}, and $n=N^{\gamma}$ for some fixed $\gamma>0$. In all cases, we use the baseline $P_i = 1/N$ $\forall i$ for $H_0$, i.e., $P$ is the uniform distribution over $N$ categories. Experiments with other baseline rates $(P_1,\ldots,P_N)$ lead to similar results provided $\min_i {n P_i}/{\log(N)} 
$ is large in the High Counts Case ($\dense$), while $\max_i \frac{n P_i}{\log(N)}$ is small in the Low Counts Case ($\sparse$). To simulate Case $(\dense)$, we use $\gamma = \xiDense$, so that $nP_i/\log(N) \approx 6.05$. To simulate Case $(\sparse)$, we use $\gamma = \xiSparse$, so that $nP_i/\log(N) \approx 0.0014$. We consider the HC and the min-P test statistics in the two-sample problem \eqref{eq:hyp}. \par
%
\subsection{Empirical Power and Phase Transition}
For each test statistic $T = T_{n,N}$ and each Monte-Carlo simulation configuration, we construct an $\alpha$-level test using as critical value $\hat{t}_{1-\alpha,M}$, the $1-\alpha$ empirical quantile of $T$ under the null hypothesis $H_0$. To determine this threshold, we simulate $M=1000$ instances under $H_0$. Next, for each configuration $(\beta,r)$, we generate $M=1000$ problem instances according to $H_1$. We define the (Monte-Carlo simulated) \emph{power} of the test statistic $T$ as the fraction of instances in which $T$ exceeds its associated threshold $\hat{t}_{1-\alpha}$. We denote this power by $\hat{B}(T,\alpha,\beta,r)$. 
\par
In order to evaluate the \emph{empirical phase transition} of the test statistic, we first indicate whether the power of $T$ is significant at each point $(\beta,r)$ in our configuration. We say that $\hat{B}(T,\alpha,\beta,r)$ is \emph{substantial} if we can reject the hypothesis 
\begin{align*}
    H_{\alpha}~~:~~ M \cdot \hat{B}(T,\alpha,\beta,r) \sim \Bin(M,\alpha).
\end{align*}
We declare $\hat{B}(T,\alpha,\beta,r)$ substantial if 
\[
\Pr \left( \Bin(M,\alpha) \geq M \cdot \hat{B}(T,\alpha,\beta,r) \right) \leq 0.05.
\]
Next, we fix $\beta \in I_\beta$ and focus on the strip $\{(\beta,r),\, r \in I_r\}$. We construct the binary-valued vector indicating those $r$ for which $\hat{B}(T,\alpha,\beta,r)$ is substantial. To this vector, we fit the logistic response model 
\[
 \Pr\left( \hat{B}(T,\alpha,\beta,r) \text{ substantial} \right) = \sigma(r|\theta_0(\beta),\theta_1(\beta)) \equiv \frac{1}{1+e^{-(\theta_1(\beta) r + \theta_0(\beta))}}.
\]
The \emph{phase transition} point of the strip $\{(\beta,r),\, r \in I_r\}$ is defined as the point $r^*(\beta)$ at which $\sigma(r^*(\beta)|\theta_0(\beta),\theta_1(\beta))=1/2$. The empirical phase transition curve is defined as $\{ r^*(\beta),\, \beta \in I_{\beta} \}$.

\subsection{Results}
\begin{figure}
\begin{center}
\input{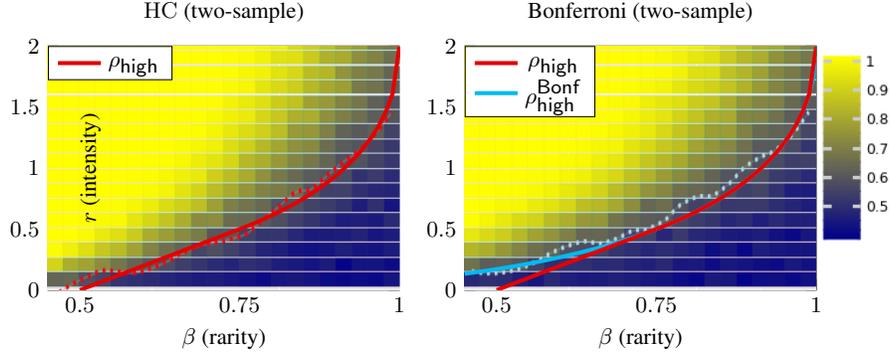}
\caption{Empirical phase diagram (High Counts Case). Shaded attribute depicts Monte-Carlo simulated power ($\Pr_{H_1}(\text{reject $H_0$})$) at the level $\Pr_{H_0}(\text{reject $H_0$}) \leq 1-\testSizeDense$ with $N=10^\npow$, $n=N^{\xiDense}$, and $P_i = 1/N$ (the uniform distribution over $N$ elements), for the $\HC$ test (left) and the min-P-value test (right). The solid red curve depicts $r=\rho_{\dense}(\beta)$, the theoretical phase-transition of HC in the two-sample setting. The dashed line represents the fitted empirical phase transition.
}
\label{fig:sim_dense}
    \end{center}
\end{figure}

\begin{figure}
\begin{center}
\input{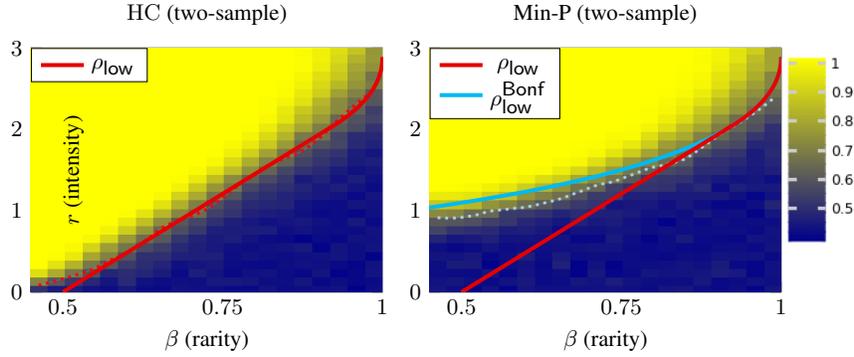}
\caption{Empirical phase diagram ($\mathsf{Sparse}$ Case). Shaded attribute depicts Monte-Carlo simulated power ($\Pr_{H_1}(\text{reject $H_0$})$) at the level $\Pr_{H_0}(\text{reject} H_0) \leq 1-\testSizeSparse$ with $N=10^\npow$, $n=N^{\xiSparse}$, and $P_i = 1/N$ (the uniform distribution over $N$ elements), for the $\HC$ (left) and min-P test (right) in the two-sample setting \eqref{eq:hyp}. The red solid curve depicts $r=\rho_{\sparse}(\beta)$, the theoretical phase transition of the $\HC$ test. The blue solid curve depicts  $r=\rho_{\sparse}^{\Bonf}(\beta)$, the theoretical phase transition of the min-P-value test. The dashed line represents the fitted empirical phase transition. 
}
\label{fig:sim_sparse}
    \end{center}
\end{figure}

Figures~\ref{fig:sim_dense} and \ref{fig:sim_sparse} illustrate the Monte-Carlo simulated power and the empirical phase transition curve in the dense and sparse cases. The results illustrated in these figures support our theoretical finding in Theorems~\ref{thm:HC_power1} and \ref{thm:HC_power2}, establishing the curves $\rho_{\dense}(\beta)$ and $\rho_{\sparse}(\beta)$ as the boundary between the region where HC has maximal power and the region where it has no power. Also shown in these figures is the Monte-Carlo simulated power and the empirical phase transition for the min-P-value test in each case. 


\section{Discussion} \label{sec:discussion}

Theorems~\ref{thm:HC_power1} and \ref{thm:HC_power2} characterize the asymptotic performance in the two-sample problem \eqref{eq:hyp}, for higher criticism of the binomial allocation P-values $\pi_1,\ldots,\pi_N$ of \eqref{eq:Pval_def}. A key enabler of our characterization is the distinction between the High Counts Case ($\dense$) and Low Counts Case ($\sparse$), as the behavior of HC varies dramatically between cases. \par
In the High Counts Case, the behavior of HC resembles its behavior in the two-sample normal means model \eqref{eq:two_sample_normal_means}, as shown by Theorem~\ref{thm:two_sample_normal_means}. Below, we conjecture that underlying this resemblance is an asymptotic equivalence between the two models. 
\par
The situation is considerably more interesting in the Low Counts Case. The phase transition curve $\rho_{\sparse}(\beta)$ appears to be new. As opposed to the High Counts Case, this curve does not seem to correspond to the phase transition of any previously known simple model. We conjecture that applying HC to the binomial allocation P-values $\pi_1,\ldots,\pi_N$, does not have the optimal phase diagram in the very sparse case where $nP_i \ll 1$. We explore this topic in a future work.

\subsection{Equivalence between Poisson and normal models?}
Our results in Theorems~\ref{thm:Bonferroni_dense} and \ref{thm:two_sample_normal_means} imply that, in the High Counts Case. the same phase transition curve $\rho_{\dense}(\beta)$ describes both the asymptotic power of HC applied to P-values deriving from two-sample normal means tests in \eqref{eq:two_sample_normal_means} and to P-values deriving from two-sample binomial tests in the Poisson means model \eqref{eq:hyp}. This equality reminds us of other results, 
showing a kind of asymptotic 
equivalence between statistical experiments in the two models. 
In the results we are thinking of, one applies a variance-stabilizing transformation $\Scal(x) \triangleq 2\sqrt{x}$ to Poisson data \cite{efron1982transformation,mccullagh1989generalized,brown2001root}, yielding data which 
are approximately normally distributed in the high-counts limit.
Applied to the data $(X_i)$ and $(Y_i)$ of \eqref{eq:hyp}, this transformation leads to the calibration
\begin{align}
\begin{split}
    \label{eq:two-sample_normal_parm}
\nu_i & = \Scal(n P_i) = 2\sqrt{n P_i}, \quad i=1,\ldots,N. \\
\nu_i^\pm & = \Scal(nQ_i^\pm) = 2\sqrt{nP_i} \pm \sqrt{2 r \log(N)},\quad i=1,\ldots,N,
\end{split}
\end{align}
for \eqref{eq:two_sample_normal_means}, which is consistent with \eqref{eq:two_sample_normal_means_perturbation}. Nussbaum and Klemel{\"a} \cite{nussbaum2006constructive} considered a transformation of the form $\Scal$ with additional randomness to show that the problem of density estimation associated with a Poisson model is asymptotically equivalent, in the sense of Le Cam \cite{le2012asymptotic}, to estimating in a Gaussian sequence model. One might 
conjecture that in the High Counts Case a similar transformation can establish an asymptotic equivalence between the Poisson experiment \eqref{eq:hyp} and the associated Gaussian experiment  \eqref{eq:two_sample_normal_means}. Assuming this were established, one might further conjecture that the equivalence holds in a sufficiently strong sense that it implies equality of phase diagrams. We have preferred here the more direct route of 
determining phase diagrams, leaving the theory of equivalent phase diagrams an intriguing area for future work.


\subsection{Randomization of P-value
}
Heuristically, it is easiest to understand the use of the HC test and its analysis if the P-values follows a uniform distribution under the null. Indeed, the Higher Criticism is well-motivated as a goodness-of-fit test of the P-values against the uniform distribution \cite{donoho2004higher,jager2007goodness}. 
For discrete situations like the Poisson means model, the P-values used in practice are stochastically larger than uniform under the null. 
Decision theory suggests to randomize the P-values of a discrete model so that their distribution is exactly uniform under the null. For the purpose of establishing a region where HC has full power as done in this paper, the decision theorist's randomized P-values and the practitioner's non-randomized P-values are asymptotically equivalent. Namely, the same phase transition emerges if we analyze HC of the randomized P-values obeying \eqref{eq:Pval_def_rand} instead of \eqref{eq:Pval_def}. 

Empirically, we observed that non-randomized P-values have some benefit in terms of power over randomized P-values. Consequently, in practice, we recommend to use the non-randomized P-values due to this empirical observation and theoretical asymptotic equivalence.

\subsection{Samples of unequal sizes}
\newtext{
It is sometimes desirable to consider the case where the two samples have unequal sizes, since in many applications we are interested in identifying chafnges in the generating mechanisms of frequency tables while ignoring differences in the total counts. In our mathematical framework, the case of unequal sample sizes corresponds to the hypothesis testing problem}
\newtext{
\begin{align}
    \begin{split}
    H_0^{(n_x,n_y)} &\quad : \quad X_i \sim \Pois(n_xP_i),\quad Y_i \sim \Pois(n_yP_i),\quad i=1,\ldots,N.\\
    H_1^{(n_x,n_y)} &\quad : \quad X_i \sim \Pois(n_xP_i) \quad \text{and} \\
    &Y_i \sim (1-\epsilon)\Pois(n_yP_i) + \frac{\epsilon}{2} \Pois(n_yQ_i^+) + \frac{\epsilon}{2} \Pois(n_yQ_i^-)
    \quad i=1,\ldots,N.
        \label{eq:hyp_unequal}
    \end{split}
\end{align}
}
\newtext{
P-values for this problem are obtained by modifying $\pi(x,y)$ of \eqref{eq:pval_def_pre} by replacing $p$ with $p'= n_x/(n_x + n_y)$. Of course, $n_x$ and $x_y$ are unknown in most applications, hence we propose to use the data-dependent estimate of $p'$,}
\begin{align*}
\hat{p}' = \frac{\hat{n}_x}{\hat{n}_x + \hat{n}_y},\quad \text{where}\quad
\hat{n}_x \equiv \sum_{i=1}^N X_i\quad \text{and} \quad \hat{n}_y \equiv \sum_{i=1}^N Y_i,
\end{align*}
\newtext{or some modifications of $\hat{p}'$ as in \cite{kipnis2019higher}.}

\newtext{
As for characterizing the asymptotic power of HC under \eqref{eq:hyp_unequal}, some situations already follow from our analysis. As an example, suppose that we are in the High Counts Case and that $n_y$ is much larger than $n_x$ that the underlying Poisson rates of each $Y_i$ is known with high accuracy. Specifically, we require that the approximation
\begin{align*}
    \pi(X_i,Y_i) & = 
    \Pr \left( \left| \Bin(X_i+Y_i, \hat{p}') - (X_i+Y_i)\hat{p}' \right| \leq \left| X_i (1-\hat{p}') - \hat{p}' Y_i \right| \right) \\
    & \approx \Pr \left( \left| \Pois(n_x P_i) - n_x P_i \right| \leq \left| X_i - n_x P_i \right| \right),
\end{align*}
is accurate on the moderate deviation scale. If we calibrate $(\epsilon, \mu)$ to $n=n_x$ in \eqref{eq:calibration}, our setting in this case reduces to the one-sample problem discussed in Section~\ref{sec:one_sample} for which HC's phase transition is given by $\rho_{\onesample}(\beta)$. More generally, consider \eqref{eq:hyp_unequal}
in the High Counts Case provided we set $n = \min\{n_x,n_y\}$ in \eqref{eq:calibration}. HC of the modified P-values is certainly asymptotically powerful whenever $r > \rho_{\twosample}(\beta)$. However, depending on the asymptotic relationship between $n_x$ and $n_y$, HC may also be asymptotically powerful for $\rho_{\onesample}(\beta) < r \leq  \rho_{\twosample}(\beta)$.}

\section{Proofs}
\label{sec:main_proof}

\subsection{Technical Lemmas}
This section provides a series of technical lemmas to be used in the proofs of our main results below.

\begin{lemma}{\cite[Lem. 3]{arias2015sparse}}\label{lem:Chernoff_bound}
Let the random variable $\Upsilon \sim \Pois(\lambda)$ and set $h(x)\equiv x \log(x)-x+1$. Then
\[
 -\lambda h(\lceil x \rceil/\lambda)-\frac{1}{2} \log \lceil x \rceil -1 \leq \log \Pr(\Upsilon_\lambda \geq x ) \leq -\lambda h(x/\lambda), \quad x>\lambda>0,
\]
and 
\[
 -\lambda h(\lceil x \rceil/\lambda)-\frac{1}{2} \log \lceil x \rceil -1 \leq \log \Pr(\Upsilon_\lambda \leq x ) \leq -\lambda h(x/\lambda), \quad 0\leq x<\lambda.
\]
\end{lemma}

Define
\begin{align} \label{eq:alpha(q,r)_def}
\alpha_{\sparse}(q,r) \equiv q\frac{ \log \left(\frac{2q}{ r\log(2)}\right)-1}{\log (2)} +\frac{r}{2}.
\end{align}

\begin{lemma} \label{lem:pi_q_sparse}
Let $\Upsilon_{\lambda}, \Upsilon_{\lambda'}$ be two independent Poisson random variables with rates $\lambda = \lambda(N)$ and $\lambda' = \lambda'(N)$, respectively. Assume that $\lambda' = \lambda + \frac{1}{2} r \log(N)(1+o(1))$, where $\lambda/\log(N) \to 0$. Fix $q>0$. Then:
\[
\Pr\left(\pi(\Upsilon_{\lambda},\Upsilon_{\lambda'}) \leq N^{-q} \right) = N^{-\alpha_{\sparse}(q,r)(1+o(1))}.
\]

\end{lemma}

\subsubsection*{Proof of Lemma~\ref{lem:pi_q_sparse}.}
From the definition of $\pi(x,y)$ in \eqref{eq:Pval_def}, we have
\begin{align*}
& \pi(\Upsilon_{\lambda},\Upsilon_{\lambda'}) = \Pr\left(\Bin(\Upsilon_{\lambda}+\Upsilon_{\lambda'},1/2) \leq \Upsilon_{\lambda} \right) 
 + \Pr\left(\Bin(\Upsilon_{\lambda}+\Upsilon_{\lambda'},1/2) \geq \Upsilon_{\lambda}' \right), 
\end{align*}
hence,
\[
\Pr(\pi(\Upsilon_{\lambda},\Upsilon_{\lambda'}) < s) \leq  
\Pr \left\{\Pr\left(\Bin(\Upsilon_{\lambda}+\Upsilon_{\lambda'},1/2) \leq \Upsilon_{\lambda} \right) < s \right\}.
\]
In addition, $\lambda<\lambda'$ implies that 
\[
\Pr\left(\Bin(\Upsilon_{\lambda}+\Upsilon_{\lambda'},1/2) \geq \Upsilon_{\lambda}' \right) < \Pr\left(\Bin(\Upsilon_{\lambda}+\Upsilon_{\lambda'},1/2) \leq \Upsilon_{\lambda} \right), 
\]
thus
\[
\Pr(\pi(\Upsilon_{\lambda},\Upsilon_{\lambda'}) < s) \geq  
\Pr \left\{\Pr\left(\Bin(\Upsilon_{\lambda}+\Upsilon_{\lambda'},1/2) \leq \Upsilon_{\lambda} \right) < s/2 \right\}.
\]
Let $y^*(x,t)$ be the threshold for $Y_i$ above which a binomial allocation P-value \eqref{eq:Pval_def} with $X_i=x$
is smaller than $t$. Namely, 
\begin{align}
    \label{eq:y_thresh}
y^*({x,t}) \equiv \arg\min_{y} \left\{y>x,\, \pi(x,y)\leq t \right\}.
\end{align}
Note that $y^*(x,t)$ is the $1-t$ quantile of the negative binomial distribution with number of failures $x$ and probability of success $1/2$. $y^*(x,t)$ is non-decreasing in $x$ and non-increasing in $t$. In addition, $y^*(0,t) = \log_2(2/t)$.
\par
Lemma~\ref{lem:Chernoff_bound} implies that 
\[
\Pr\left( \Upsilon_{\lambda'} \geq y^*(0,t) \right) \leq \exp\{-\lambda'h(y^*(0,t)/\lambda')\}.
\]
From $y^*(0,t) = \log_2(2/t)$, and $h(p)=p\log p-p+1$,
\[
\lambda' h(y^*(0,N^{-q})/\lambda') = q \log_2(2N) \left(\log \frac{q \log_2(2N)}{\lambda'} - 1 \right) + \lambda',
\]
where we set $t = N^{-q}$. Now, as $N\to \infty$, 
\[
\frac{\lambda'}{\log(N)} \sim \frac{r}{2},
\]
\[
\log \left( \frac{q \log_2(2N)}{\lambda'} \right) \to \log \left( \frac{2q}{r \log(2)}\right).
\]
It follows that
\[
\frac{q \log(2N)}{\log(N)} \left( \log \frac{q \log(2N)}{\lambda'}-1\right) + \frac{\lambda'}{\log(N)} \sim \frac{q}{\log(2)}\left( \log \frac{2q}{r\log(2)} -1 \right) + \frac{r}{2},
\]
and 
\[
\frac{\lambda' h\left(y^*(0,N^{-q})/\lambda' \right)}{\log(N)} \sim \alpha_{\sparse}(q,r),\quad N\to \infty.
\]
From here,
\begin{align*}
  \Pr\left(\pi(\Upsilon_{\lambda},\Upsilon_{\lambda'}) \leq N^{-q} \right) & \leq 
 \sum_{x=0}^\infty \Pr\left( \Upsilon_{\lambda} = x \right) \left( \Pr \left( \Upsilon_{\lambda'} \geq y^*(x,N^{-q}) \right) \right) \\
 & \leq
 \sum_{x=0}^\infty \Pr\left( \Upsilon_{\lambda} = x \right) \Pr \left( \Upsilon_{\lambda'} \geq y^*(0,N^{-q}) \right) \\
 & = \Pr \left( \Upsilon_{\lambda'} \geq y^*(0,N^{-q}) \right) = N^{-\alpha_{\sparse}(q,r)+o(1)}.
\end{align*}

\par We now develop a lower bound on $\Pr\left(\pi(\Upsilon_{\lambda},\Upsilon_{\lambda'}) \leq N^{-q} \right)$, starting from an upper bound on $y^*(x,t)$. We have
\begin{align*}
     \pi(x,y) & \leq 2\left(2^{-(x+y)} \sum_{k=0}^x \binom{x+y}{k} \right) \\
    & \leq 2^{-(x+y)+1} (1+x+y)^x \leq 2^{-y} (2+4y)^x,
\end{align*}
where the last transition follows from $x<y$ for $N^{-q}<1/2$, valid as $N\to \infty$. 
The condition
\[
2^{-y} (2+4y)^{x} \leq t,
\]
implies that $2^{-y^*} (2+4y^*)^{x} \leq t$ where $y^* = y^*(x,t)$. Assuming that $x \leq  \lceil \lambda \rceil /a_N$, $a_N = \log \log(N)$, $\lambda/\log(N)\to 0$, $t=N^{-q}$, and $ \log(N) \leq y$, we get
\[
y^*(x,t) \leq q\log_2(N)(1+o(1)).
\]
From Lemma~\ref{lem:Chernoff_bound} in the case $x < \lambda$, and using $\lambda/\log(N)\to 0$, we obtain:
\begin{align*}
    \Pr \left( \Upsilon_{\lambda} \leq \lceil \lambda \rceil/a_N \right) &= N^{-o(1)}.
\end{align*}
We use the above lower bounds on $y^*(x,t)$ and $\Pr \left( \Upsilon_{\lambda} \leq \lceil \lambda \rceil/a_N \right)$ in the following:
\begin{align*}
 & \Pr\left(\pi(\Upsilon_{\lambda},\Upsilon_{\lambda'}) \leq N^{-q} \right)
 \geq  \Pr\left\{ \Pr\left(\Bin(\Upsilon_{\lambda},\Upsilon_{\lambda'},1/2) \leq \Upsilon_{\lambda}  \right)  \leq N^{-q}/2 \right\} \\
 & \qquad = 
 \sum_{x=0}^\infty \Pr\left( \Upsilon_{\lambda} = x \right) \Pr \left( \Upsilon_{\lambda'} \geq y^*(x,N^{-q}) \right) \\
  & \qquad \geq \sum_{x \leq \lceil \lambda \rceil} \Pr\left( \Upsilon_{\lambda} = x \right) \Pr \left( \Upsilon_{\lambda'} \geq y^*(x,N^{-q}) \right) \\
 & \qquad \geq \Pr(\Upsilon_{\lambda} \leq \lambda)
\Pr \left(\Upsilon_{\lambda'} \geq
q\log_2(N)(1+o(1))
\right) \\
 & \qquad \overset{(a)}{\geq}  N^{-o(1)} \exp\left(-q \log_2(N)(1+o(1)) \left(\log\frac{q \log_2(N)(1+o(1))}{\lambda'} -  1\right)  \right. \\
 & \qquad \qquad \left. + \lambda' - \frac{1}{2} \log \lceil q  \log_2(N) (1+o(1)) \rceil  - 1 \right) \\
& \qquad  = N^{-o(1)} \exp\left(-q \log_2(N)(1+o(1)) \left(\log\frac{q \log_2(N)(1+o(1))}{\lambda'} -  1\right) \right) \\
 & \qquad = N^{-\alpha(q,r)(1+o(1))}.
\end{align*}
where $(a)$ follows from Lemma~\ref{lem:Chernoff_bound} and $\Pr(\Upsilon_\lambda \leq \lambda) >1/3$. \qed

\begin{lemma} \label{lem:modorate_deviation}
Let $\Upsilon_\lambda',\Upsilon_\lambda$ denote two independent Poisson random variables. Let $a(\cdot)$ be a real-valued function, $a(x) : (0,\infty) \to (0,\infty)$. Consider a sequence of pairs $(\lambda, \lambda') = (\lambda_N, \lambda'_N)$ where each $\lambda'_N > \lambda_N$. Suppressing subscript $N$, suppose 
$\lambda \to \infty$, $\lambda' \geq \lambda$, $\lambda'/\lambda \to 1$. Also suppose $a(\lambda) - (\sqrt{2\lambda'} - \sqrt{2\lambda}) \to \infty$ while $a(\lambda)/\lambda \to 0$. Then:
\[
\lim_{\lambda \to \infty} \frac{1}{ \left(\sqrt{a(\lambda)}
-(\sqrt{2\lambda'}-\sqrt{2\lambda})
\right)^2}  \log \left[ \Pr\left( \sqrt{2\Upsilon_{\lambda'}}-\sqrt{2\Upsilon_{\lambda}} \geq \sqrt{a(\lambda)} \right) \right] = -\frac{1}{2}.
\]
\end{lemma}
\subsection{Proof of Lemma~\ref{lem:modorate_deviation}}
By normal approximation to the Poisson, as $\lambda \to \infty$, 
\[
\frac{\Upsilon_{\lambda}-\lambda}{\sqrt{\lambda}} \overset{D}{\to} \Ncal(0,1).
\]
The transformed random variable $\sqrt{\Upsilon_\lambda}$ is asymptotically variance-stabilized \cite{anscombe1948transformation,mccullagh1989generalized}: 
\[
2(\sqrt{\Upsilon_{\lambda}}-\sqrt{\lambda}) \overset{D}{\to} \Ncal(0,1),  \quad \lambda \to \infty;
\]
Because $\log(n)/\lambda \to 0$, our result is a consequence of a ``moderate deviation'' estimate (see \cite{RubinSethuraman1965} \cite[Ch. 3.7]{zeitouni1998large}) for the random variable
\[
\sqrt{2\Upsilon_{\lambda'}}- \sqrt{2\Upsilon_{\lambda}} + \sqrt{2\lambda'} - (\sqrt{2\lambda'} + \sqrt{2\lambda}).
\]
\qed


\begin{lemma} \label{lem:yx_lower}
Consider a sequence $\{\lambda_N\}$ such that $\lambda_N/\log(N) \to \infty$. Fix $q>0$, and define $\tilde{y} \equiv \tilde{y}_{N,q} \equiv \left( \sqrt{x} + \sqrt{ q\log(N)-a_N} \right)^2$,  
where $\{a_N\}$ denotes a positive 
sequence satisfying $a_N \lambda \geq \log^2(N)$ and $a_N/\log(N) \to 0$. There exists $N_0(q)$ so that for all $N > N_0(q)$, and $x\geq \lambda - \sqrt{a_N \lambda}$,
\[
\pi(x,\tilde{y}) \geq N^{-q}.
\]

\end{lemma}
\subsubsection*{Proof of Lemma~\ref{lem:yx_lower}}
As $x,y\to \infty$ in such a way that $y/x \to 1$, we have
\begin{align*}
    &\pi(x,y) = 2^{-(x+y)-1} \sum_{k=0}^x \binom{x+y}{k} \geq 2^{-(x+y)-1} \binom{x+y}{x}\\
    & \qquad  = 2^{-(x+y)-1} \frac{(1+o(1))}{\sqrt{2 \pi}} \sqrt{\frac{x+y}{xy}} \left(1+\frac{y}{x} \right)^{x} \left(1+\frac{x}{y} \right)^{y};
\end{align*}
the last step by Stirling's approximation. 
Set $x^* = \lambda - \sqrt{a_N \lambda}$. Since $x\leq \tilde{y}$, we get
\begin{align*}
& \inf_{x \geq \lambda - \sqrt{a_N \lambda}} N^q \pi(x, \tilde{y}) = N^{q} \pi(x^*,\tilde{y}) \\
& \qquad  \geq N^q 2^{-(x^*+\tilde{y})-1} \frac{(1+o(1))}{\sqrt{2 \pi}} \sqrt{\frac{x^*+\tilde{y}}{x^*\tilde{y}}} \left(1+\frac{\tilde{y}}{x^*} \right)^{x^*} \left(1+\frac{x^*}{\tilde{y}} \right)^{\tilde{y}}.
\end{align*}
The proof is completed by verifying that, under our assumptions on $\{a_N\}$, the last expression goes to infinity as $\lambda$ goes to infinity. \qed

\begin{lemma} \label{lem:pi_q_dense}
Let $\Upsilon_{\lambda}, \Upsilon_{\lambda'}$ be two independent Poisson random variables with rates $\lambda$ and $\lambda'$, respectively. Assume that $\lambda' = \lambda + \sqrt{2 \lambda r \log(N)(1+o(1))}$ and  $\lambda/\log(N) \to \infty$. Fix $q>r/2$. Let $\alpha(q,r) \equiv \alpha_{\dense}(q,r) = (\sqrt{q}-\sqrt{r/2})^2$. Then,
\[
\Pr\left(\pi(\Upsilon_{\lambda},\Upsilon_{\lambda'}) \leq N^{-q} \right) \geq N^{-\alpha(q,r)+o(1)}
\]
and 
\[
\Pr\left(\pi(\Upsilon_{\lambda},\Upsilon_{\lambda'}) \leq N^{-q} \right) \leq 
N^{-1+o(1)} + N^{-\alpha(q,r)+o(1)}.
\]

\end{lemma}

\subsubsection*{Proof of Lemma~\ref{lem:pi_q_dense}} 
Consider the threshold level $y^*(x,t)$ of \eqref{eq:y_thresh}. 
Hoeffding's inequality \cite{hoeffding1994probability}
\[
\Pr( \Bin(n,p) \leq t) \leq e^{\frac{-2(t-pn)^2}{n}}, 
\]
implies
\[
\pi(x,y) = 2\Pr \left( \Bin(x+y,1/2) \leq x \right) \leq 2 e^{-\frac{(y-x)^2}{2(x+y)}}
\]
for all integers $y\geq x \geq 0$. 
Therefore, the conditions
\[
\frac{(y-x)^2}{x+y} \geq 2\log(2/t), \quad y>x,
\]
imply $\pi(x,y) \leq t$. Because $y^*(x,t) \geq x$ for  $0<t<1/2$, we solve for $y$ and get 
\begin{align*}
y^*(x,t) & \leq x + \log(2/t) + 2\sqrt{x \log(2/t) + (\log(2/t))^2} 
\end{align*} 
whenever $t<1/2$. Set $t=N^{-q}$ for some fixed $q>0$. We have
\begin{align*}
 \Pr & \left(\pi(\Upsilon_{\lambda},\Upsilon_{\lambda'}) \leq N^{-q} \right) \\
 & \geq 
 \Pr \left( \Upsilon_{\lambda'} \geq 
 \Upsilon_{\lambda} + q\log(2^{1/q}N) + 2\sqrt{q\Upsilon_{\lambda} \log(2^{1/q}N) + (q\log(2^{1/q}N))^2} 
 \right) \\
 & = 
 \Pr \left( \Upsilon_{\lambda'} \geq 
 \Upsilon_{\lambda} + q\log(N)(1+o(1)) + 2\sqrt{q\Upsilon_{\lambda} \log(N)(1+o(1))}
 \right) \\
 & \overset{(a)}{=} \Pr \left( \sqrt{\Upsilon_{\lambda'}}-\sqrt{\Upsilon_{\lambda}} \geq \sqrt{q\log(N)(1+o(1))} \right)  \\
 & = \Pr \left( \sqrt{2\Upsilon_{\lambda'}}-\sqrt{2\Upsilon_{\lambda}} \geq \sqrt{2q\log(N)(1+o(1))} \right),
\end{align*}
where (a) follows from the equivalence of the events $\sqrt{\Upsilon_\lambda'}\geq \sqrt{\Upsilon_\lambda}+\sqrt{\Delta}$ with $\Upsilon_{\lambda'} \geq \Upsilon_\lambda + 2 \sqrt{\Upsilon_{\lambda} \Delta} + \Delta$, where $\Delta \equiv q \log(N)$. Because $\sqrt{2\lambda'}-\sqrt{2\lambda} = \sqrt{r\log(N)(1+o(1))}$,
Lemma~\ref{lem:modorate_deviation} implies that, for each fixed $q>r/2$,
\[
\Pr \left( \sqrt{2\Upsilon_{\lambda'}}-\sqrt{2\Upsilon_{\lambda}} \geq \sqrt{2q\log(N)(1+o(1))} \right) = N^{-(\sqrt{q}-\sqrt{r/2})^2 + o(1)}.
\]

\par For the upper bound, use Lemma~\ref{lem:yx_lower} to conclude that the threshold level \eqref{eq:y_thresh} satisfies
\begin{align*}
& y^*(x,N^{-q}) \geq  (\sqrt{x} + \sqrt{q \log(N)(1+o(1))})^2
\end{align*}
for all $x$ such that $x \geq \lambda - \sqrt{a_N \lambda}$, where $\{a_N\}$ satisfies $a_N \lambda \geq \log^2(N)$. We obtain
\begin{align*}
 \Pr & \left(\pi(\Upsilon_{\lambda},\Upsilon_{\lambda'}) \leq N^{-q} \right) =
\Pr \left( \Upsilon_{\lambda'} \geq y^*(\Upsilon_{\lambda},N^{-q}) \right) \\
& = \Pr \left( \Upsilon_{\lambda'} \geq y^*(\Upsilon_{\lambda},N^{-q}) \mid \Upsilon_{\lambda} \geq \lambda- \sqrt{a_N \lambda} \right) \Pr \left(\Upsilon_{\lambda} \geq \lambda- \sqrt{a_N \lambda} \right) \\
& \qquad +
\Pr \left( \Upsilon_{\lambda'} \geq y^*(\Upsilon_{\lambda},N^{-q}) \mid \Upsilon_{\lambda} < \lambda- \sqrt{a_N \lambda} \right) \Pr \left( \Upsilon_{\lambda} < \lambda- \sqrt{a_N \lambda} \right) \\
 & \leq \Pr \left( \sqrt{\Upsilon_{\lambda'}}-\sqrt{\Upsilon_{\lambda}} \geq \sqrt{q \log(N)(1+o(1))} \right) +  \Pr \left(\Upsilon_{\lambda} < \lambda - \sqrt{a_N \lambda} \right),
\end{align*}
where we have used the equivalence of the event $\Upsilon_\lambda' \geq (\sqrt{\Upsilon_\lambda} + \sqrt{q \log(N)(1+o(1))})^2$ with $\sqrt{\Upsilon_\lambda'} - \sqrt{\Upsilon} \geq \sqrt{q \log(N)(1+o(1))}$ to get 
\[
\Pr \left( \Upsilon_\lambda' \geq y^*(\Upsilon_\lambda,N^{-q}) \mid \Upsilon_\lambda \geq \lambda - \sqrt{a_N \lambda} \right) \leq \Pr \left( \sqrt{\Upsilon_{\lambda'}}-\sqrt{\Upsilon_\lambda} \geq \sqrt{q \log(N)(1+o(1))} \right). 
\]
Now, Lemma~\ref{lem:Chernoff_bound} leads to
\begin{align*}
& \log \Pr\left(\Upsilon_{\lambda} < \lambda - \sqrt{a_N \lambda}\right) \leq -(\lambda - \sqrt{a_N \lambda})\log\left({1-\sqrt{\frac{a_N}{\lambda}}} \right) + \sqrt{a_N \lambda}  \\
& = -\sqrt{a_N \lambda}(1+o(1)) \leq -\log(N)\cdot(1+o(1)).
\end{align*}
Lemma~\ref{lem:modorate_deviation} implies
\begin{align*}
    & \Pr \left( \sqrt{\Upsilon_{\lambda'}}-\sqrt{\Upsilon_{\lambda}} \geq \sqrt{q \log(N)(1+o(1))} \right) =
     N^{-(\sqrt{q}-\sqrt{r/2})^2+o(1)}.
\end{align*}
It follows that
\begin{align*}
 & \Pr\left(\pi(\Upsilon_{\lambda},\Upsilon'_{\lambda'}) \leq N^{-q} \right) \leq N^{-1+o(1)} + N^{-(\sqrt{q}-\sqrt{r/2})^2+o(1)}.
\end{align*}
\qed

The following lemma characterizes the behavior of $\HC_{N,n}^*$ under $H_0$, by comparing it to the normalized uniform empirical process. 
\begin{lemma} \label{lem:under_null}
Under $H_0$ of \eqref{eq:hyp}, we have
\[
\Pr \left( \HC_{N,n}^\star \leq \sqrt{4 \log \log (N)} \right) \to 1.
\]
\end{lemma}
\subsubsection*{Proof of Lemma~\ref{lem:under_null}}
Let $U_1,\ldots,U_N$ be i.i.d. samples from the uniform distribution on $(0,1)$. Denote by 
\begin{align*}
& F_N(t) \equiv \frac{1}{N} \sum_{i=1}^N \one\{\pi_i \leq t\}, \quad F^{(0)}_N(t) \equiv \frac{1}{N} \sum_{i=1}^N \one \{U_i \leq t\},
\end{align*}
the empirical distribution of $\pi_1,\ldots,\pi_N$ and $U_1,\ldots,U_N$, respectively. The normalized uniform empirical process
\[
W_N(t) \equiv \sqrt{N} \frac{  F^{(0)}_N(t)-t }{\sqrt{t(1-t)}},
\]
is known to satisfy \cite{shorack2009empirical}
\[
 \frac{\max_{0<t \leq \alpha_0} W_N(t)}{\sqrt{2 \log \log N}} \overset{p}{\to} 1, 
\]
as $N\to \infty$. We have
$\Pr_{H_0}(\pi_i \leq t) \leq t$ and $\Pr(U_i \leq t) = t$, hence $F_N(t) \leq F^{(0)}_N(t)$ stochastically.
Since
\[
\HC_{N,n}^* = \max_{0<t \leq \alpha_0} \sqrt{N}\frac{F_N(t)-t}{\sqrt{t(1-t)}}, 
\]
we get that, stochastically, 
\[
\HC_{N,n}^* \leq \max_{0<t \leq \alpha_0} W_N(t),
\]
and hence
\[
\Pr_{H_0} \left( \HC_N^* \leq \sqrt{4 \log \log (N)} \right) \to 1.
\]
\qed

The following lemma provides an asymptotic lower bound on the empirical cumulative distribution function of the P-values $\pi_1,\ldots,\pi_N$ under $H_1$ for values tending to zero on the scale $N^{-q}$, $q\in(0,1)$. 
\begin{lemma} \label{lem:main_lemma}
Let $\alpha(\cdot)$ and $\gamma(\cdot)$ be two real-valued functions $\alpha, \gamma : [0,\infty) \to [0,\infty)$. Let $q\in(0,1)$ and $\beta>0$ be fixed. Suppose that $\alpha(q)<q$, and that $F_{N,n}$ satisfies
\begin{align}
    \label{eq:main_lemma}
\ex{F_{N,n}(N^{-q})} = N^{-q+o(1)}(1-N^{-\beta}) +
N^{-\beta} N^{-\alpha(q)+o(1)}.
\end{align}
Let $\{a_N\}$ be a positive sequence obeying $a_N N^{-\eta} \to 0$ for any $\eta>0$. 
    If 
    \begin{align}
    \label{eq:main_lem_assump_1}
    \alpha(q)+\beta < \gamma(q),
\end{align}
then 
\[
\Pr( N^{\gamma(q)} (F_{N,n}(N^{-q}) -N^{-q}) \leq a_N) = o(1).
\]
\end{lemma}
\subsubsection*{Proof of Lemma~\ref{lem:main_lemma}}
Set $t_N=N^{-q}$ and $\eta = \gamma(q) - \alpha(q)-\beta>0$. We have
\begin{align*}
    & \Pr \left( F_{N,n}(t_N) - t_N \leq a_N N^{-\gamma(q)} \right) \\
    & \qquad = \Pr \left( F_{N,n}(t_N) - t_N \leq (1-\delta) (\ex{F_{N,n}(t_N) - t_N}) \right),
\end{align*}
where $\delta = \delta_N$ obeys:
\begin{align}
\delta & = 1 - \frac{a_N N^{-\gamma(q)}}{\ex{F_{N,n}(t_N)-t_N}} \nonumber \\
    & = 1 - \frac{ a_N N^{-\gamma(q)}}{  
    N^{-q}(N^{o(1)}-N^{-\beta})
    + N^{-\beta - \alpha(q)+o(1)}
    }
    \nonumber \\
    & = 1 - \frac{ a_N N^{-\eta}}{  
    N^{\alpha(q)-q}(N^{\beta+o(1)}-1)
    + N^{o(1)}
    }. \label{eq:lem:proof1}
\end{align}
The assumptions $\alpha(q)>q$ and $a_N N^{-\eta}\to 0$ imply that, eventually, $0<\delta<1$. For $X$ the sum of $N$ independent Bernoulli random variables with $\mu = \ex{X}$, the Chernoff inequality \cite[Ch 4.]{mitzenmacher2017probability} says
\[
\Pr \left( X \leq (1-\delta)\mu \right) \leq \left( \frac{e^{-\delta}}{(1-\delta)^{1-\delta}} \right)^\mu \leq e^{-\mu \frac{\delta^2}{2}},\qquad \delta\in (0,1].
\]
We use this inequality with $X = N F_{N,n}(t) = \sum_{i=1}^N \one \{\pi_i \leq t\}$, which leads to
\begin{align*}
    & \log \Pr \left( F_{N,n}(t_N) - t_N \leq a_N N^{-\gamma(q)} \right) \leq  -\frac{\delta^2 N}{2} \ex{F_{N,n}(t)-t_N} \\
    & \qquad  \overset{(a)}{=} -\frac{\delta N}{2}  \left(\ex{F_{N,n}(t)-t_N} - a_N N^{-\gamma(q)} \right) \\
    & \qquad \overset{(b)}{=} -\frac{\delta}{2} \left( N^{1-q}( N^{o(1)}-N^{-\beta+o(1)}) + N^{1-\beta-\alpha(q)} (N^{o(1)} - a_N N^{-\eta}) \right). 
\end{align*}
where (a) follows by \eqref{eq:lem:proof1} and we invoked \eqref{eq:main_lemma} in step (b). 
Since $q<1$ and $a_N N^{-\eta} \to 0$, this last expression goes to $-\infty$ for any fixed choice of $\alpha$ and $\beta$. 
\qed

\subsection{Region of No Power \label{sec:powerlessness}}
Consider the randomized P-values $\tilde{\pi}_1,\ldots \tilde{\pi}_n$ of \eqref{eq:Pval_def_rand}. We have 
\begin{align}
\tilde{H}_0~~:~~ \tilde{\pi}_i \simiid \Unif(0,1),\quad i=1,\ldots,n,
\end{align}
and
\begin{align}
\label{eq:H1}
\tilde{H}_1~~:~~ \tilde{\pi}_i \simiid (1-\epsilon_N)\Unif(0,1) + \epsilon_N G_{N,i},\quad i=1,\ldots,n,
\end{align}
where $G_{N,i}$ is a continuous probability distribution such that, for $X_i \sim G_{N,i}$, 
\begin{align*}
\max_{i=1,\ldots,N}\left( -\log \left(\Pr(X_i < N^{-q})\right) \right) \leq \log(N) \left( \alpha_{*}(q,r) + o(1) \right), \qquad * \in \{\dense, \sparse\}.
\end{align*}
Under these conditions, powerlessness of HC for $r < \rho_*(\beta)$ F following result of \cite{DonohoKipnis2020b}.
\begin{theorem}{\cite[Thm. 2.1]{DonohoKipnis2020b}} \label{thm:converse}
Consider testing $\tilde{H}_0$ against $\tilde{H}_1$. Suppose that $G_{N,i}$ has a continuous density $f_i$, and that, for some $C< \infty$, $\|f_i\|_{\infty} = C$ for all $i=1,\ldots,N$, and for $X_i \sim G_{N,i}$ 
\begin{align}
    \label{eq:tail}
\max_{i=1,\ldots,N}\left( -\log \left(\Pr(X_i < N^{-q})\right) \right) \leq \log(N) \left( \alpha(q,r) + o(1) \right), 
\end{align}
for some bivariate function  $\alpha(q,r)$ that is continuous, non-negative, increasing in $q$ and decreasing in $r$. $\HC_N^*$ is asymptotically powerless if \[
\max_{0\leq q \leq 1} \left(\frac{q+1}{2} - \alpha(q,r) - \beta \right) < 0. 
\]
\end{theorem}

\subsection{Region of Full Power}
\label{sec:full_power}

Over the domain $q \in (0,1+\eta)$, $\eta>0$, and $r>0$, suppose that $\alpha_{*}(q,r)$ is a bivariate continuous non-negative function, non-increasing in $r$ and non-decreasing in $q$. Let $\rho_{*}(\beta)$ be the minimal $r$ satisfying
\[
\max_{0\leq q \leq 1} \left( \frac{q+1}{2} - \beta - \alpha_{*}(q,r) \right) \leq 0,
\]
and hence $r > \rho_{*}(\beta)$ implies that there exists $q = q_0(\beta,r) \in (0,1)$ so that 
$(q_0+1)/2 > \beta + \alpha_{*}(q,r)$. In the rest of the proof, let $q = q_0(\beta,r)$. \par
Under $H_1$, Lemma~\ref{lem:pi_q_dense} implies that \begin{align*}
 & \Pr_{H_1} \left( \pi_i \leq N^{-q}\right)= (1-\epsilon_N)N^{-q+o(1)} + \epsilon_N \Pr(\pi(\Upsilon_{nP_i},\Upsilon_{nQ_i^+}) \leq N^{-q}) \\
 & =  (1-\epsilon_N)N^{-q+o(1)} + \epsilon_N N^{-\alpha_{*}(q,r)+o(1)},
\end{align*}
uniformly in $i$, where $\alpha_{*}(q,r) = \left(\sqrt{q}-\sqrt{r/2}\right)^2$ (in the last display and throughout the proof, $o(1)$ represents an expression satisfying $\max_i o(1)\to 0$ as $N\to \infty$). It follows that 
\[
\ex{F_{N,n}(N^{-q})} = (1-\epsilon_N)N^{-q+o(1)} + \epsilon_N N^{-\alpha_{*}(q,r)+o(1)}.
\]
In view of Lemma~\ref{lem:under_null}, it is enough to show that, as $N \to \infty$,
\begin{align} \label{eq:to_show}
\Pr_{H_1} \left( \HC_{N,n}^\star \leq \sqrt{4\log \log(N)} \right) \to 0,
\end{align}
whenever $r > \rho_{*}(\beta)$. 
Using $t_N=N^{-q}$ and $a_N = \sqrt{4 \log\log(N)}$,
\begin{align}
    \Pr \left( \HC_{N,n}^\star \leq a_N \right) & \leq \Pr \left( \sqrt{N} \frac{F_{N,n}(t_N)-t_N}{\sqrt{t_N (1-t_N)}} \leq a_N \right) \nonumber  \\
    & \leq \Pr \left( N^{\frac{q+1}{2}}(F_{N,n}(t_N) -t_N) \leq a_N \right). \label{eq:proof_main1}
\end{align}
Next, apply Lemma~\ref{lem:main_lemma}
to \eqref{eq:proof_main1} with $\alpha(q) = \alpha_{*}(q,r)$, $\gamma(q) = (q+1)/2$, and $a_N = \sqrt{4 \log \log(N)}$. This lemma yields \eqref{eq:to_show}. 

\subsection{Proof of Theorem~\ref{thm:HC_power1}}
Substituting $* = \dense$, 
Part (i) follows from Section~\ref{sec:full_power} and Part (ii) follows from Section~\ref{sec:powerlessness}. \qed

\subsection{Proof of Theorem~\ref{thm:HC_power2}}
Substituting $* = \sparse$, 
Part (i) follows from Section~\ref{sec:full_power} and Part (ii) follows from Section~\ref{sec:powerlessness}. \qed


\subsubsection*{Proof of Theorem~\ref{thm:Bonferroni_dense}}
First note that the condition $r > \rho_{\dense}^\Bonf(\beta)$ is equivalent to \begin{equation} \label{eq:Bonf_r_cond}
1 > (1- \sqrt{r/2})^2 + \beta.
\end{equation}
For any $t\in(0,1)$, we have:
\begin{align*}
\Pr_{H_0}(\pi_{(1)} \leq t) \leq 1-(1-t)^N.
\end{align*}
Pick a sequence $\{t_N\}$ obeying $n t_N \to 0$. Along this sequence, the last display goes to zero. Below we use the specific sequence $t_N =
(2N \log(N))^{-1}$. Let $\Upsilon_{\lambda_1}$ and $\Upsilon_{\lambda_2}$ denote two independent Poisson random variables with rates $\lambda_1$ and $\lambda_2$, respectively. We have
\begin{align*}
    & \log \Pr_{H_1} \left( \pi_{(1)} > t_N \right) = \sum_{i=1}^N \log \Pr_{H_1}\left(\pi_i > t_N \right),
\end{align*}
and 
\begin{align*}
    \Pr_{H_1}\left(\pi_i > t_N \right) & = (1-\epsilon_N) \Pr_{H_0}\left(\pi_i > t_N \right) + \frac{\epsilon_N}{2} \Pr \left( \pi(\Upsilon_{nP_i},\Upsilon_{nQ_i^+}) > t_N \right)  \\
    & \qquad + \frac{\epsilon_N}{2} \Pr \left( \pi(\Upsilon_{nP_i},\Upsilon_{nQ_i^-}) > t_N \right) \\
    & = (1-\epsilon_N) \Pr_{H_0}\left(\pi_i > t_N \right) + \epsilon_N \Pr \left( \pi(\Upsilon_{nP_i},\Upsilon_{nQ_i^+}) > t_N \right) \\
    & \leq 1-\epsilon_N + \epsilon_N\left(1 - \Pr\left( \pi(\Upsilon_{nP_i},\Upsilon_{nQ_i^+}) \leq t_N \right) \right).
\end{align*}
To complete the proof, it is enough to show that, for $r$ and $\beta$ satisfying \eqref{eq:Bonf_r_cond}, for $t_N \equiv (2N \log(N))^{-1}$, and for every $i=1,\ldots,N$,
\[
N \epsilon_N  \Pr(\pi(\Upsilon_{nP_i},\Upsilon_{nQ_i^+}) \leq t_N) \to \infty, 
\]
since this would imply that $\Pr_{H_1} \left( \pi_{(1)} > t_N \right) \to 1$. 
Let $y^*({x,t})$ be the threshold level \eqref{eq:y_thresh}. 
For Poisson random variables $\Upsilon, \Upsilon'$, 
\[
\Pr \left( \pi (\Upsilon',\Upsilon) > t \right)= \sum_{x=0}^\infty \Pr(\Upsilon = x)\Pr\left( \Upsilon' \geq y^*(x,t) \right).
\]
It follows that 
\begin{align*}
    & \epsilon_N  \Pr(\pi(\Upsilon_{nP_i},\Upsilon'_{nQ_i^+}) > t_N) =
 \epsilon_N  \sum_{x=0}^\infty  \Pr(\Upsilon_{nP_i}=x)
  \Pr(\Upsilon'_{nQ_i^+} \geq y^*(x,t_N)) \\
  & \geq \epsilon_N 
    \sum_{x \leq nP_i+\sqrt{nP_i}} \Pr(\Upsilon_{nP_i}=x) \Pr\left(\Upsilon'_{nQ_i^+} \geq y^*(nP_i+\sqrt{nP_i},t_N) \right),
\end{align*}
where we used the fact that $y^*(x,t)$ is non-decreasing in $x$ for $t<1/2$. Using the same Chernoff bound argument as in the proof of Lemma~\ref{lem:pi_q_dense}, we also have 
\begin{align*}
y^*(x,t) & \leq x + \log(2/t) + 2\sqrt{x \log(2/t) + (\log(2/t))^2} 
\end{align*} 
whenever $t<1/2$. Consequently, 
\begin{align*}
 \Pr &\left(\pi(\Upsilon_{nP_i},\Upsilon'_{nQ_i^+}) \leq t_N \right) = \Pr\left(\Upsilon'_{nQ_i^+} \geq y^*(\Upsilon_{nP_i}, t_N) \right)
 \\
 & \geq 
 \Pr \left( \Upsilon'_{nQ_i^+} \geq 
 \Upsilon_{nP_i} + \log(2/t_N) + 2\sqrt{ \Upsilon_{nP_i} \log(2/t_N) + (\log(2/t_N))^2} \right) \\
 & =
 \Pr \left( \Upsilon'_{nQ_i^+} \geq 
 \left(\sqrt{\Upsilon_{nP_i}} + \sqrt{\log(2/t_N)(1 + o(1))} \right)^2 \right) \\
 & = \Pr \left( \Upsilon'_{nQ_i^+} \geq 
 \left(\sqrt{\Upsilon_{nP_i}} + \sqrt{\log(N)(1 + o(1))} \right)^2 \right).
\end{align*}
From
 \[
\sqrt{\log(2/t_N)} - (\sqrt{2 nQ_i^+} - \sqrt{2 nP_i}) = \sqrt{\log(N)+\log(\log(N))}- \sqrt{r \log(N)} \to \infty,
\]
Lemma~\ref{lem:modorate_deviation} implies 
\begin{align*}
 \log \Pr & \left( \Upsilon_{nQ_i^+} \geq 
 \left(\sqrt{\Upsilon_{nP_i}} + \sqrt{\log(2/t_N)(1 + o(1))} \right)^2 \right)  \\
 & \qquad = -\frac{1}{2} \left(\sqrt{2\log(N)\log(\log(N))}- \sqrt{r \log(N)} \right)^2 + o(1) \\
 & \qquad = -\log(N) \left( (1-\sqrt{r/2})^2 + o(1) \right).
\end{align*}
Hence, from \eqref{eq:Bonf_r_cond},
\[
N \epsilon_N  \Pr\left(\pi(\Upsilon_{nP_i},\Upsilon_{nQ_i^+}) \leq t_N \right) \geq N^{1-\beta-(1-\sqrt{r/2})^2+o(1)} \to \infty.
\]
\qed

\subsubsection*{Proof of Theorem~\ref{thm:Bonferroni_sparse}}
Arguing as in the proof of Theorem~\ref{thm:Bonferroni_dense}, it is enough to show that, for $Nt_N \to 0$, and for every $i=1,\ldots,N$,
\[
N \epsilon_N  \Pr(\pi(\Upsilon_{nP_i},\Upsilon_{nQ_i^+}) \leq t_N) \to \infty. \]
Consider the threshold level $y^*(x,t)$ of \eqref{eq:y_thresh}. As in the proof of Lemma~\ref{lem:pi_q_sparse}, from
\begin{align*}
    & \pi(x,y) \leq  2^{-y} (2+4y)^x,
\end{align*}
and $t_N= 1/(N \log(N))$, we get that 
\[
y^*(x,(N \log(N))^{-1}) \leq  \log_2(N)(1+o(1)),
\]
whenever $ \log(N) \leq y$ and  $x\leq \lceil 2n P_i \rceil$. 
We have
\begin{align*}
    \Pr &(\pi(\Upsilon_{nP_i},\Upsilon_{nQ_i^+}) \leq t_N) = \sum_{x=0}^\infty \Pr(\Upsilon_{nP_i}=x) \Pr(\Upsilon_{nQ_i^+} \geq y^*(x,t_N)) \\
    & \geq 
    \sum_{x\leq \lceil 2nP_i \rceil } \Pr(\Upsilon_{nP_i}=x) \Pr(\Upsilon_{nQ_i^+} \geq \log_2(N)(1+o(1))) \\
    & \overset{(a)}{\geq} \frac{1}{2} \exp\left(- \log_2(N)(1+o(1)) \left(\log\frac{ \log_2(N)(1+o(1))}{\lambda'} -  1\right)  \right. \\
&\qquad \left. + \lambda' - \frac{1}{2} \log \lceil  \log_2(N) (1+o(1)) \rceil  - 1 \right) \\
& = \frac{1}{2} \exp\left(- \log_2(N)(1+o(1)) \left(\log\frac{ \log_2(N)(1+o(1))}{\lambda'} -  1\right) \right) = N^{-\alpha_{\sparse}(1,r)+o(1)},
\end{align*}
where $(a)$ follows from Lemma~\ref{lem:Chernoff_bound} and from $\Pr(\Upsilon_\lambda \leq 2 \lambda ) \geq 1/2$. We conclude that 
\[
N \epsilon_N \Pr(\pi(\Upsilon_{nP_i},\Upsilon_{nQ_i^+}) \leq t_N) \to 0
\]
whenever $1-\beta - \alpha_{\sparse}(1,r) > 0$, which is equivalent to the condition $ r >  \rho_{\sparse}^{\Bonf}(\beta)$.  \qed

\subsection{Proof of Theorem~\ref{thm:two_sample_normal_means}}

As in the case of Theorems~\ref{thm:HC_power1} and \ref{thm:HC_power2}, the key to characterizing the power behavior of Higher Criticism is a lemma on behavior of the individual P-values. Theorem~\ref{thm:two_sample_normal_means} is based on the following result:
\begin{lemma}
\label{lem:normal_Pvals}
Let $X \sim \Ncal(\nu,1)$ and $Y\sim \Ncal(\nu',1)$ be independent. Set
\[
\bar{\pi}(x,y) \equiv 
\Pr\left( \left|\Ncal(0,1)\right| \geq \frac{\left|y - x\right|}{\sqrt{2}} \right),
\]
and assume that $\nu' = \nu \pm \sqrt{2r\log(N)}$ and $q > r/2$. Then:
\[
\Pr\left( \bar{\pi}(X,Y) \leq N^{-q} \right) = N^{-(\sqrt{q}-\sqrt{r/2})^2+o(1)}.
\]
\end{lemma}

\subsubsection*{Proof of Lemma \ref{lem:normal_Pvals}}
Set $U = (Y-X)/\sqrt{2}$ and note that $U\sim N(\sqrt{r \log(N)},1)$.
Standard facts about Mills' ratio imply
\begin{align}
    \label{eq:Mill's_ratio}
\Pr \left(\left|\Ncal(0,1)\right| \geq |x| \right) \sim \frac{2\phi(x)}{|x|} = e^{-\frac{x^2}{2}(1+o(1))},\quad x \to \infty.
\end{align}
Therefore, for $Z\sim \Ncal(0,1)$ and $\mu_N \to \infty$ as $N\to \infty$, 
\begin{align*}
&  \Pr \left(\left|\Ncal(0,1)\right| \geq Z + \mu_N \right) = e^{-\frac{(Z+\mu_N)^2}{2}(1+o_p(1))}. 
\end{align*}
We have
\begin{align*}
 \Pr( \bar{\pi}(X,Y) \leq N^{-q}) & = \Pr  \left(\Pr \left(\left|\Ncal(0,1)\right| \geq U \right) \leq N^{-q} \right) \\
  & = \Pr  \left(\Pr \left(\left|\Ncal(0,1)\right| \geq Z + \sqrt{ r \log(N)} \right) \leq N^{-q} \right) \\
& = \Pr \left( e^{-\frac{(Z+\sqrt{r \log(N)})^2}{2}(1+o_p(1))} \leq e^{-q \log(N)} \right) \\
& = \Pr \left( Z \geq \log(N) \left(\sqrt{2q}-\sqrt{r} \right)(1 + o_p(1)) \right) \\
& = N^{-(\sqrt{q}-\sqrt{r/2})^2+o(1)},
\end{align*}
where in the last transition we used $2q>r$ and \eqref{eq:Mill's_ratio}.  Applying Lemma~\ref{lem:normal_Pvals} to the P-values $\bar{\pi}_1,\ldots,\bar{\pi}_N$ of \eqref{eq:normal_Pvalues}, we see that 
\[
\Pr_{H_1}( \bar{\pi}_i \leq N^{-q}) = (1-\epsilon_N)N^{-q+o(1)} + \epsilon_N N^{-\alpha_{\dense}(q,r) + o(1)}. 
\]
From here, the proof is identical to the proof of Theorem~\ref{thm:HC_power1}. \qed


\bibliographystyle{imsart-number}  
\bibliography{HigherCriticism}   

\end{document}